\documentclass[preprint,12pt,authoryear]{elsarticle}

\usepackage{amssymb}
\usepackage{amsmath}

\usepackage{lineno}

\usepackage{booktabs}
\usepackage{float}
\usepackage[hidelinks]{hyperref}
\usepackage{cleveref}
\usepackage{listings}
\usepackage{multirow}

\journal{Journal of Construction Engineering and Management}

\begin{document}

\begin{frontmatter}



\title{Advancing Strategic Planning and Dynamic Control of Complex Projects} 


\author[a]{L.G. Teuber}
\author[b]{H.J. van Heukelum}
\author[c]{A.R.M. Wolfert\corref{cor1}}\ead{a.r.m.wolfert@tudelft.nl}

\cortext[cor1]{Corresponding author}

\affiliation[a]{organization={Data Science \& Engineering, Boskalis},
            addressline={Rosmolenweg 20}, 
            city={Papendrecht},
            postcode={3356 LK}, 
            state={Zuid-Holland},
            country={The Netherlands}}

\affiliation[b]{organization={Corporate R\&D, Boskalis},
            addressline={Rosmolenweg 20}, 
            city={Papendrecht},
            postcode={3356 LK}, 
            state={Zuid-Holland},
            country={The Netherlands}}

\affiliation[c]{organization={Department Algorithmics, Faculty of Mathematics and Computer Science, Delft University of Technology},
            addressline={Van Mourik Broekmanweg 6}, 
            city={Delft},
            postcode={2628 XE}, 
            state={Zuid-Holland},
            country={The Netherlands}}

\begin{abstract}
\noindent Strategic project planning and dynamic control are essential to ensure that complex projects are both prepared and executed best-fit-for-common-purpose, guided by three interrelated strategies: (1) Agreeing First, (2) Acting Feasibly, and (3) Adapting Flexibly.  When these strategies become too complex for humans to fully conceive and manage, effective computer-aided decision support becomes crucial. However, standard simulation-driven evaluation and \textit{a-posteriori} decision-making are typically single-sided and technically focused focus, rather than applying a combined simulation-and-optimisation approach that \textit{a-priori} integrates stakeholder interests and their mitigation behavior. Moreover, current planning and control methodologies often lack robust stochastic representations and associative multi-objective optimisation methods that capture the full socio-technical complexity while maximizing the potential within reach. This paper introduces Odycon (Open Design and Dynamic Control), a new purpose-driven project management methodology that provides an actionable solution to these challenges. It presents a generic mathematical framework for project planning and control that integrates stakeholder preferences (human domain) with system performances (physical domain), enabling more effective planning and dynamic control. Odycon integrates standard Monte Carlo simulation (MCS) with the novel Integrative Maximisation of Aggregated Preferences (IMAP) optimisation method to develop a best-fit strategic plan and the most effective mitigation control strategies. Its value is demonstrated through applications in offshore wind installation and highway infrastructure projects, showcasing advances in associative design and decision-making, and aiming for a best-fit-for-common-purpose synthesis across different complex project phases.
\end{abstract}




\end{frontmatter}

\section{Introduction}

\noindent In today's complex engineering project management landscape, wicked problems often cause projects to derail. In fact, projects don't go wrong — they start wrong \citep{Flyvbjerg2024}. And even when they start right, they can still lose direction. Advanced design and project management practices may still hit impasses because project managers often rely on retrospective analyses - such as those by \citet{Flyvbjerg2024} or others like \citet{Yang2022} - and traditional decision-support methods, which offer little to no actionable guidance for a constructive, best-for-project way forward. These challenges become particularly acute when decision-makers fail to recognize that their problems are part of a larger, interconnected whole, and that mutual interconnection and collaboration are key to unlocking complex projects that achieve best-fit for a common purpose.

This calls for moving beyond the traditional notion of \emph{fitness for purpose} — often associated with total quality bound by scope, cost, and time (also known as project management triangle or triple constraint) and the delivery of technically effective results that meet predefined requirements. In today’s complex construction environments, this narrow, end-product-oriented view is no longer sufficient. Instead, we must adopt a broader, value-driven perspective that focuses on what is socially desirable and physically feasible guided by a common purpose from the very start and continuing throughout to the end. \emph{Best-fit for common purpose} integrates stakeholder interests, dynamic context, and collective value creation within socio-physical reach.

Moreover, what is needed is a proactive threefold approach to confront project complexity — ensuring that projects start, proceed, and finish successfully: (1)  \textit{Agreeing First} on a best-fit scope, (2) \textit{Acting Feasibly} with an optimal project plan, and (3) \textit{Adapting Flexibly} through the best set of mitigation measures. Establishing a robust and socio-technically feasible project plan is fundamental to effective and efficient project management \citep{Smith2014}. However, as project size and uncertainty increase, initial estimates regarding scope, cost, and timelines are often exposed to inherent risks and uncertainties that cannot be entirely eliminated.

Project success is therefore achieved through dynamic adaptation to emerging challenges—often not specified in static contracts—and through collaborative engagement among stakeholders or decision-makers to co-develop associative solutions. Focusing solely on the final scope as defined in the contract, insisting on creating a fixed plan and adhering rigidly to it, is ultimately naïve and pseudo-realistic in the face of real-world project dynamics.

This ultimately calls for an actionable, open design and systems-oriented approach, supported by transparent decision-support tools that foster both deliberate (slow) and socially informed thinking. As such, successful project delivery hinges on purpose-driven and optimal strategic planning and dynamic control, transforming complexity into outcomes that best fit within socio-physical reach. \emph{Planning} refers to the process of reaching agreement first in order to act feasibly, while project \emph{control} ensures execution according to plan and the ability to adapt flexibly when needed on the run.

Despite this reality, traditional project management in construction often remains heavily focused on cost and schedule, with stakeholders primarily considering these from individual perspectives. This narrow focus neglects the interplay of stakeholder-oriented, concurrent objectives. Joint project success can only be enhanced when stakeholders are willing to look beyond self-interest and engage in a process of mutual alignment. This requires stakeholders to let go of their individually defined objectives in favour of a synthesis solution that best fits the common purpose — one that ultimately outperforms single-sided or even compromise-based project outcomes \citep{Wolfert2023, Scharmer2016, Glasl1998, VandenDoel1993}. 

By integrating diverse stakeholder interests and accounting for project capabilities in the decision-making process, a more holistic outcome—referred to as the best-fit for common purpose—can be achieved, rather than settling for suboptimal compromises \citep{Zhilyaev2022}. This best-fit design emerges from the socio-technical interplay between stakeholder preferences (the human domain—what they \textit{want}) and system performance capabilities (the physical domain—what it \textit{can} deliver), as developed and demonstrated by \citet{VanHeukelum2024}.

However, achieving such a best-fit requires more than technical alignment; it also demands a fundamental shift in how planning and execution control are approached. Both academic research and industry practice tend to rely on retrospective, analysis-oriented methods (\textit{a-posteriori}), which offer limited support for forward-looking strategic development. Moreover, these methods are often driven by individual optimization or simulation-based approaches, producing one-sided solutions that overlook the multifaceted nature of real-world projects. In practice, technical capabilities, human goal-oriented behavior, and stakeholder desirability are deeply interconnected. Yet traditional planning and control approaches fail to explicitly capture the dynamic interplay of all stakeholder interests within this complex decision-making process \citep{Blanchard2011, DelPico2023, slack2010operations}.

Most existing design optimisation methods rely on single-objective optimisation (SOO), often reducing multiple project objectives to a single cost metric. This approach overlooks the qualitative nature of stakeholder preferences, which are essential for goal-oriented decision-making. While the benefits of multi-objective optimisation (MOO) in project management are increasingly recognized \citep{Guo2022}, current methodologies rarely integrate MOO within strategic project planning or dynamic project control frameworks that also incorporate probabilistic risk management \citep{Kammouh2022}. To address the rising complexity and uncertainty in construction projects, decision support must evolve to enable computer-aided, stakeholder-oriented planning and control based on MOO.  As \citet{kahneman2011thinking} emphasizes, deliberate and data-informed processes are essential for mitigating cognitive biases such as overconfidence. Consequently, well-informed project decision-making requires transparent, evidence-based tools that support “slow and social thinking”—tools that not only process data but also help align diverse stakeholder objectives within feasible and resilient project solutions.

In this paper, all the aforementioned considerations are taken into account in the development of a new open-design-driven dynamic project planning and control methodology, named Odycon (acronym: Open Design and Dynamic Control). Odycon builds on the core principles and design optimization method of the Open Design Systems (Odesys) methodology \citep{Wolfert2023}. The Odesys integrative maximisation of aggregated preferences (IMAP) optimization method has proven effective in achieving best-fit-for-common-purpose design solutions, where all stakeholder interests are accounted for and translated into a common preference domain to identify an overall group optimum. Notably, Odesys/IMAP is a novel preference-based optimization methodology where, instead of directly minimizing or maximizing individual objectives, the rigorous mathematical optimization focuses solely on maximizing the aggregated preference of stakeholder objectives within physical reach while removing fundamental modelling errors of standard multi-objective optimization methods \citep{VanHeukelum2024}.

So far, the Odesys/IMAP method enables iterative group design and decision-making in a deterministic manner. However, for a \emph{project planning},  it must be extended to incorporate traditional probabilistic project planning techniques that capture schedule performance, such as PERT or discrete-event simulation (DES). This constitutes a truly novel project management-by-design methodology.  Additionally, these tools should employ stochastic simulation methods to reflect the inherent uncertainties and unpredictability encountered both in the planning and execution phases. For \emph{project control} applications he method must also be expanded to shift the focus from optimizing design variables to optimizing the allocation of control measures, using multiple project objectives to evaluate their effectiveness and efficiency. To achieve this, the state-of-the-art project manager-oriented mitigation controller (MitC), as developed and demonstrated by \citet{Kammouh2022, Kammouh2021}, will be enhanced from a single-cost objective orientation to optimally facilitate stochastic, multi-objective project control. In doing so, Odycon constitutes a truly novel \emph{human-centred project management-by-design} methodology. 

To summarize, the following Odycon advances are presented in this paper:

\begin{enumerate}
    \item A general framework for integrated strategic planning and dynamic control, combining probabilistic project schedule simulations with preference-based multi-objective design optimisation. This enables dynamic adaptation to project-specific conditions and stakeholder objectives, thereby enhancing associative decision-making in complex construction environments.
    \item A direct linkage between project feasibility, project management objectives, and stakeholder decisions, resulting in a unified integration of network performance, stakeholder values, and associative planning and control processes.
    \item A human-centred, actionable project management decision-support model that overcomes fundamental modelling errors in classical multi-objective optimisation and Monte Carlo planning by fully integrating human preferences with physical performances, eliminating invalid aggregation and Pareto-front assumptions, and accounting for human-oriented mitigation behavior \citep{Kammouh2021, VanHeukelum2024}. 
    \item A pure preference-based aggregation mechanism for multi-objective decision-making in project management. This translates diverse objectives and constraints into a shared preference domain, enabling an optimised group synthesis without relying on monetisation \citep{Wolfert2023}. Rather than selecting the lowest cost solution, the approach supports identifying the solution with the highest combined value for both project outcomes and people—that is, a best-fit for common purpose plan and set of control measures.
    \item A methodological advancement from PII to PIII in project management practice. Building on the PII model for complex construction projects introduced by \citet{VanGunsteren2011}, itself rooted in double-loop learning as described by \citet{Argyris1996}, this paper proposes an evolution toward PIII—Best Practice for Collaborative Projects. Unlike PII, which focuses on making improved decisions for internal stakeholders via integrative reflection and learning, PIII promotes associative open-loop learning and co-creation with all relevant stakeholders. This fosters dynamic group-based planning and control, where integrative maximisation of associated preferences replaces static, self-oriented decision-making. In doing so, it enables a genuinely collaborative approach to achieving the best possible overall outcome for both project and people \citep{Wolfert2023}.
\end{enumerate}

Combining the advances of MitC and IMAP enables a comprehensive representation of all stakeholder preferences and objectives in project planning and control, while simultaneously accounting for socio-physical constraints. This significantly reduces the risk of inefficient or misinformed decision-making—especially as project size and complexity increase, making it virtually impossible to manually consider all potential scenarios. Odycon supports the generation of best-fit-for-common-purpose solutions for both strategic planning and dynamic control applications. As a decision-support methodology, Odycon automates the selection of planning and control variables by integrating multiple stakeholder objectives and constraints into a unified stochastic simulation and optimisation process.\\

This paper is structured as follows. It begins by introducing the general mathematical formulation of Odycon’s multi-objective project planning and control optimisation approach. Next, the underlying framework and core concepts of the stochastic simulation and optimisation methodology are explained, including a detailed description of the method’s design, within which the IMAP optimisation is embedded. The methodology is first demonstrated using a strategic planning case - \textit{Agreeing First \& Acting Feasibly } (with no control variables) - involving an offshore wind installation project, where a service provider and a marine contractor are the key decision-makers. This is followed by a dynamic control case - \textit{Adapting Flexibly }(with no planning variables) - applied to a highway infrastructure project. This demonstrator features associative mitigation control by both a contractor and a highway agency, illustrating how to identify a best set of control measures to enable flexible adaptation across two distinct project management phases. The paper concludes with recommendations for further development and reflects on the added value of the newly developed stochastic simulation and optimisation methodology.

\section{Mathematical formulation of the multi-objective optimisation}\label{sec:optimisation}

\noindent To address the mentioned shortcomings of project planning \& control, a new stochastic optimisation \& simulation methodology is established that enables the intricate representation of all stakeholders interests towards the different project objectives. Reflecting the goal-oriented behaviour of involved stakeholders is achieved by integrating the human preference domain (associated stakeholder interests) and the physical planning performance domain using the Odesys methodology embedded within a stochastic simulation framework \citep{Wolfert2023,VanHeukelum2024}. In other words, design-oriented stakeholder behaviour needs to be integrated into probabilistic planning performance. To this end, the integrative multi-objective design optimisation method IMAP, as part of the Odesys methodology, must first be formulated for a multi-objective planning \& control optimisation problem. This threefold integration of performance-, objective- and preference functions, reads as follows:

\begin{equation}
    \label{eq:general_MS}
    \begin{gathered}
    \mathop{Maximise}_{\mathbf{x}} \textrm{ } U = A \left[  
    P_{k,i} 
    \left(O_i
    \left(\mathcal{N}(\mathbf{x}, \mathbf{y}), F_1(\mathbf{x}, \mathbf{y}),F_2(\mathbf{x}, \mathbf{y}),...,F_J(\mathbf{x}, \mathbf{y}) \right)
    \right),w'_{k,i} \right] \textrm{ for } \\ 
    k=1,2,...,K \\ 
    i=1,2,...,I
    \end{gathered}
\end{equation}

\noindent

Subject to: 
\begin{equation}
    \label{eq:ineq_cons}
    \begin{gathered}
    g_{p}(O_{i}(\mathcal{N}(\mathbf{x}, \mathbf{y}), F_{1,2,...,J}(\mathbf{x}, \mathbf{y})),F_{1,2,...,J}(\mathbf{x}, \mathbf{y})) \le 0
    \textrm{ for } p=1,2,...,P
    \end{gathered}
\end{equation}

\begin{equation}
    \label{eq:eq_cons}
    \begin{gathered}
    h_{q}(O_{i}(\mathcal{N}(\mathbf{x}, \mathbf{y}), F_{1,2,...,J}(\mathbf{x}, \mathbf{y})),F_{1,2,...,J}(\mathbf{x}, \mathbf{y})) = 0
    \textrm{ for } q=1,2,...,Q
    \end{gathered}
\end{equation}

\noindent

\noindent
With:

\begin{itemize}
    \item $U$: Utility function that needs to be maximized for a best-fit configuration, is addressed using a Genetic Algorithm software tool, \texttt{preferendus} \citep{VanHeukelum2024}.
    
    \item $A$: Aggregated preference score determined as part of the IMAP optimisation method (see \ref{a-fine-a}).

    \item $P_{k,i}(O_{i}(\mathcal{N}(\mathbf{x}, \mathbf{y}), F_{1,2,...,J}(\mathbf{x}, \mathbf{y}))$: Preference functions that describe the preference stakeholder $k$ has towards objective function $i$, which are functions of different performance functions and dependent on planning and/or control variables and (physical) variables ( $i \le k$ and $K$ is the maximum number of stakeholders).

    \item $O_{i}(\mathcal{N}(\mathbf{x}, \mathbf{y}), F_{1,2,...,J}(\mathbf{x}, \mathbf{y}))$: Objective functions that describe the objective $i$, functions of different performance functions, planning and/or control variables and physical variables.

    \item $\mathcal{N}(\mathbf{x}, \mathbf{y})$: Performance function that describes the project network, depending on planning and/or control variables $\mathbf{x}$ and (physical) variables $\mathbf{y}$. Here $\mathcal{N}$ is a logical representation of a project planning network, composed of nodes (activities) and edges (logical links) with their respective properties.

    \item $F_{1,2,...,J}(\mathbf{x}, \mathbf{y})$:  Performance functions that describe the object behaviour, depending on the planning and/or control variables $\mathbf{x}$ and (constant) physical variables $\mathbf{y}$.

    \item $\mathbf{x}$: A vector containing the available planning \& control variables (i.e. endogenous design variables) $x_1,x_2,...,x_N$. These variables are bounded such that $lb_n \le x_n \le ub_n$, where $lb_n$ is the lower bound, $ub_n$ is the upper bound.

    \item $\mathbf{y}$: A vector containing the physical variables (i.e. exogenous 'non-design' variables) $y_1,y_2,...,y_M$.

    \item $g_{p}(O_{i}(\mathcal{N}(\mathbf{x}, \mathbf{y}), F_{1,2,...,J}(\mathbf{x}, \mathbf{y})),F_{1,2,...,J}(\mathbf{x}, \mathbf{y}))$: Inequality constraint functions, which can be either objective function and/or performance function constraints.
    \item $h_{q}(O_{i}(\mathcal{N}(\mathbf{x}, \mathbf{y}), F_{1,2,...,J}(\mathbf{x}, \mathbf{y})),F_{1,2,...,J}(\mathbf{x}, \mathbf{y}))$: Equality constraint functions, which can be either objective function and/or performance function constraints.

    \item $w'_{k,i}$: Weights for each of the preference functions (sum of the number of preference functions equals $K*I$). The global weights for the relative importance of stakeholders is defined as $w_{k}$. The local weight stakeholder $k$ gives to objective $i$ is defined as $w_{k,i}$. The following formula holds: $w'_{k,i}=w_{k}\cdot w_{k,i}$, given that $\sum w'_{k,i}=\sum w_{k,i}=\sum w_k=1$. Note that in case of equivalent stakeholder decision-making \( w_k = 1/K \).
\end{itemize}

\noindent
To further elaborate on this mathematical formulation, several important remarks are made which are discussed below.

\subsubsection*{Remark 1: Preference aggregation}
\noindent The preference aggregation is performed using the principles of the PFM theory \citep{Barzilai2022}. It is an integral part of the Integrative Preference Aggregation Method (IMAP), see \citet {Wolfert2023}. A best-fit solution is identified as a feasible plan that yields the highest aggregated preference. Note that individual preferences reflect the level of interest in specific objectives. To arrive at an aggregated preference for individual objectives, the IMAP method employs the mathematical operator \textit{A}, which computes the aggregated preference score using an algorithm called the A-fine Aggregator (see \ref{a-fine-a}).

\subsubsection*{Remark 2: Preference functions}
\noindent The preference functions describe the relation of a stakeholders satisfaction towards a certain objective on a defined scale. In other words, preference — the only property of relevance in decision theory and equivalent to value — provides goal orientation to the objective function. Note: objective functions, in themselves, merely describe an objective and do not inherently imply a minimisation or maximisation direction — a distinction that sets IMAP apart from all classical multi-objective optimisation methods. The elicitation of the preference functions and associated individual weights is essential to reflect the relationship of the stakeholders interest towards a certain objective. As an initial estimation for both the preference functions input and local and global the weights the Conjoint Analysis (CA) method can be used. For more information regarding the process of preference elicitation, the reader is referred to \citet{Arkesteijn2017}. As mentioned in \citet{VanHeukelum2024} objective $O_i$ can be associated with multiple stakeholder preference functions $P_{k,i}$ (as $k \geq i$). Therefore, it is not required that a stakeholder expresses a preference for all objectives.

\subsubsection*{Remark 3: Network}
\noindent Optimisation requires representing the network planning as a Directed Acyclic Graph (DAG), which is a structure composed of nodes (project activities) and directed edges (activity interdependencies: i.e, logical links) that do not form any cycles, ensuring all paths from the initial node to any other node are open-ended and loop-free.

\subsubsection*{Remark 4: Planning \& control variables}
\noindent In the context of a pure project planning optimisation (no control measures), the variables $x_n$ are defined as stated above, allowing to take any form (e.g. integer, float, binary, etc.). However, for a pure project control optimisation (no planning variables), the variables $\mathbf{x}$ are redefined as $\mathbf{x} = \mathbf{a} \mathbf{c_i}$, that is the product of the allocation of a control variables vector $\mathbf{a}$ (which takes binary value $a_n \in \{ 0,1 \}$, where 0 reflects not allocated, or 1 reflects allocated) and the control impact vector $\mathbf{c_i}$ towards objective $i$, containing the impact of control variables $c_1,c_2,...,c_N$ directly related towards control variables $x_1,x_2,...,x_N$. Note that $x_{1..N} \in \{ 0,c_{1..N} \}$, which implies that the solution space is bounded.

\section{Odycon's stochastic simulation \& optimisation methodology}\label{methodology}
\noindent This section describes the conceptual functioning of the Odycon methodology, with its underlying stochastic simulation \& optimisation framework combining probabilistic MCS and IMAP multi-objective optimisation to enable associative project control for best-fit for common purpose project execution. The essence is to account for the inherent uncertainties of a project combined with the goal-oriented objectives of all stakeholders. The following two sections elaborate on the integrated MCS and the utilised IMAP optimisation method. 

\subsection{Probabilistic MCS approach}
\noindent To accurately reflect uncertainties and unpredictability encountered during project execution, the core of the Odycon methodology is build upon a probabilistic MCS frequently used for modelling uncertainties and capturing the stochastic behaviour of variables in project scheduling \citep{DelPico2023}. The duration uncertainty of activities, impact of risks, impact of other external influences, or control measures is quantified by the Beta-PERT distribution. This distribution is widely recognised for modelling uncertainties in project scheduling. Accordingly, the probability density functions for the varying parameters are denoted as $f(z_i; a_i, m_i, b_i)$, where $a_i$, $m_i$, and $b_i$ represent the optimistic, most likely, and pessimistic time estimates, respectively.



Within an MCS iteration, random values $z_i$ will be sampled according to the given respective probability distributions. The occurrence of a risk is modelled by a binary variable $X$, where a risk occurs $X=1$ with probability $p_e$, and a risk does not occur $X=0$ with probability $1-p_e$. In deterministic project scheduling (using the critical path method), time estimates are fixed, resulting in one static project network realisation. In contrast, with MCS used in probabilistic project planning (PERT), the structure of the project network $\mathcal{N}$ changes with each iteration based on random sampling, creating various realisations of project networks. This variability reflects the inherent uncertainty and stochastic nature modelled in the MCS, allowing for a comprehensive analysis of frequency distribution of optimal project planning and control strategies. This approach also aids in identifying the most probable project paths, critical activities, thereby enhancing decision-making and risk management strategies in project planning and execution.

\begin{figure}[H]
    \centering
    \includegraphics[width=0.75\linewidth]{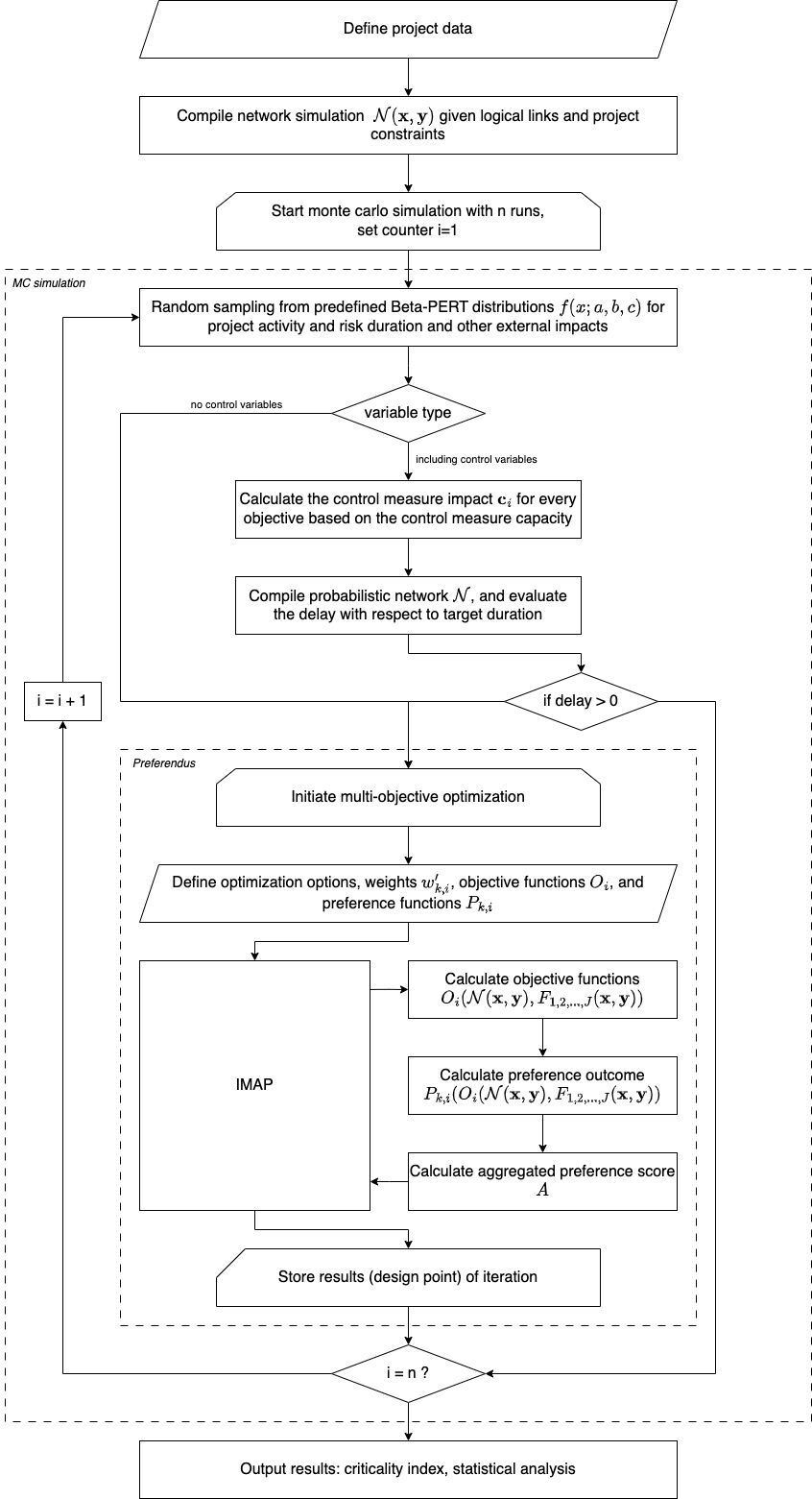}
    \caption{Odycon concept diagram.}
    \label{fig:Flowchart}
\end{figure}

The framework enables decision-making through the MCS given the steps provided in \autoref{fig:Flowchart}. (1) First, the Project network data is read, and the network $\mathcal{N}(\mathbf{x}, \mathbf{y})$ is established given the project activities, risk event and any other relevant data. Note, to that extent any network planning or simulation method can be integrated, that is able to reflect the necessary logical links in the network. (2) Next, the MCS is initiated by a loop of n runs. (3) Step 3 is to sample random duration's for the project activities, risk and other external influences (like weather). To establish the optimisation, first the type of variables must be defined. If control variables are not considered, the planning variables are defined and the optimisation initiated (continues with step 6). (4) If control variables are included, the impact $\mathbf{c_i}$ towards objective $i$ needs is calculated given the duration capacity of the control measure. (5) Next, the network is compiled and the delay with respect to the target duration evaluated. In case of a delay, the optimisation is initiated. (6) The optimisation is run by defining the objective functions, the preference functions (which defines the desired outcomes for the objective functions) and their related weights, in accordance with the prior established mathematical formulation for the optimisation (see \autoref{sec:optimisation}). Once defined, the optimisation is performed utilising the IMAP method. To arrive at the highest aggregated preference score $A$, the values of the objectives functions are calculated and with that the corresponding value of the preference functions. Once established, the aggregated preference score can be determined according to the principles of references aggregation with preference maximisation (more information see the section on the optimisation below). (7) After every iteration of the MCS the result is stored. (8) Once the MCS reaches the $n^{th}$ iteration, the results of the individual optimisations can be presented with the criticality index (how often does a certain variable or a combination of variables occur). Given the collected data, detailed statistical analysis are enabled.

\subsection{Integrative multi-objective optimisation approach}
\noindent To integrate the human goal-oriented behaviour in project management, and associative decision-making the multi-objective optimisation formulation (defined in \autoref{sec:optimisation}) is integrated into the MCS as part of the IMAP optimisation method. The optimisation is executed within every iteration of the probabilistic MCS, aiming to select the most effective set of project planning variables and or control measures given the specific realisation of the network within that iteration. The 

The essence of the IMAP method lies in combining preferences aggregation (reflecting associative decision-making) with preference maximisation (reflecting human goal-oriented behaviour). To arrive at the highest aggregated preferences the weighted (or relative) scores of the individual preference functions ($P_{k,i}$) are aggregated in their affine space according to the basic principles of PFM theory. The final preference score aggregation is performed using the weighted least squares method. The preference score aggregation as part of the IMAP optimisation is described in \ref{a-fine-a}. To facilitate the application of the IMAP method in practice, the Python-based software tool, \texttt{preferendus}, is incorporated into the framework \citep{VanHeukelum2024}.

With utilising the IMAP optimisation method the shortcoming of classical MOO methods that consider the Pareto front as a valid outcome is addressed. MOO methods using the Pareto front result in a large set of multiple possible solution points, with the need to evaluate these solution points \textit{a posteriori}. Moreover, given the complexity of a project the set of planning \& control variables and objectives can become very large, consider a for example a project with 4 stakeholders and 40 control variables. In such a case, project managers and stakeholder are interested in a single best design solution. We want to avoid cumbersome \textit{a posteriori} evaluations of possible suboptimal solutions, resulting from a trial-and-error approach. Classical Pareto front based methods are in this planning \& control approach case inappropriate (see \citet{VanHeukelum2024}). IMAP instead results in an integrative \textit{a priori} optimisation with a single best-fit result from a multidimensional solution space.

\section{Demonstrative planning \& control applications}

\noindent The introduced Odycon methodology is demonstrated in this section upon two real-life multi-stakeholder construction projects. To demonstrate the Odycon advances towards concurrent and associative design and decision-making, two examples are presented: (1) a pure strategic planning application (no control variables) for an offshore wind installation project, and (2) a pure dynamic control application (no planning variables) of a highway infrastructure construction project.

To legitimise the choice of these two examples, the following is considered. In general, two types of project management approaches are defined. For projects governed by strict procedural continuity with limited flexibility to deviate from planned tasks (the so-called Toothpaste model), management tends to prioritise strategic planning in advance. In contrast, for projects that allow considerable flexibility and adjustments during implementation (the so-called Viking model), management strategies allow primarily for dynamic control on-the-run \citep{reschke2013dimensions}.

Unlike inland highway projects where extra assets are often readily available and can be mobilised during project execution, offshore projects require detailed front-end planning due to the limited possibility of on-the-run control. For this reason, a pure planning scenario for this 'Toothpaste' project is considered where the focus is to select the most impactful strategies up-front, reducing the potential risk during execution. To show the effect of dynamic control on-the-run, a pure control scenario for a 'Viking' project is chosen, which is demonstrated on the basis of a transport infrastructure project. 
These demonstrative applications illustrate how the introduced Odycon methodology can support optimal decision-making for execution and mitigation strategies. The complexity of these projects demands computer aided design and decision-making, where the participating stakeholders together decide to arrive at best-fit for common project purpose.  

To this end, two Python based models were developed: one for a purely strategic planning case (SYPL, an acronym of strategic planning) and another for a purely dynamic control case (MICO, an acronym of mitigation control). The source code, project data, and results of these demonstrative applications are available on GitHub, as detailed in the data availability statement.

Final note, both application examples are based upon previous publications, all of which one of the authors of this paper was involved in, see \citet{VanHeukelum2024} and/or \citet{Kammouh2022}. To read these examples here independently, the project descriptions are partly repeated.

\subsection{Strategic Planning application : Agree First \& Act Feasibly}

\noindent \textit{Technical context:} Offshore floating wind (OFW) technologies are considered a solution for wind energy production in deep waters and regions with stronger wind velocities \citep{spring2020global}. Contrary to the conventional method of employing fixed monopiles for bottom-founded installations, these turbines are mounted on platforms that are anchored to the seabed using mooring lines. The OFW site consists of 36 floating wind turbines (FWT) and 108 suction anchors (i.e., 3 anchors per FWT). OFW projects are executed by a set of highly specialised vessels, where the installation of the individual foundations or floaters is often performed in an iterative or cyclic process (serial installation). The cyclic installation process is often constrained by a singular, inflexible execution path, leaving little room for deviation in the face of delays from adverse weather, technical issues, or operational risks. This rigidity necessitates direct resolution of problems, pausing the project if needed, rather than employing alternative strategies. \\

\noindent\textit{Social context:} This application includes the following concurrent stakeholder objectives both from an energy service provider (Stakeholder 1) and a marine contractor (Stakeholder 2): (1) minimising project duration, (2) reducing installation costs, (3) optimising fleet utilisation, and (4) lowering CO2 emissions. While the energy service provider prioritises fast project completion to expedite revenue generation and aims to decrease CO2 emissions for better environmental and societal acceptance, the marine contractor focuses on cutting costs to enhance competitiveness and seeks to improve fleet utilisation for greater operational efficiency.

We will now first describe the integrative planning optimisation problem by working through the mathematical statement (see \autoref{sec:optimisation}), resulting in performance-, objective-, and preference functions.

\subsubsection{Performance functions}

\noindent \textit{Planning variables}\\
The installation of OFW turbines requires specialised types of installation vessel. The performance of the planning sequence is dependent on the utilisation of these vessels during construction. The type of vessels and their properties used in the project are defined in \autoref{tab:owf_variables}. The boundaries are the (integer) number of possible vessel to be utilised. The utilisation refers to the probability of a vessel being better utilised in a different project of the contractor.

\begin{table}[H]
\caption{Strategic planning application: Planning variables (OCV refers to offshore construction vessel).}
\resizebox{\textwidth}{!}{\begin{tabular}{l c c c c c c}
\toprule
Variables & Description & Boundaries & Anchor deck space & Day rate $R_i$ & utilisation prob. $p_i$ & CO\textsubscript{2} emissions $E_i$ \\
\midrule
$F1 = x1$ &   small OCV    &   $0 \le x1 \le 3$   &   8 &  $47 k\text{\texteuro} $ & 0.7 & $30 \text{\ }t/day$ \\
$F2 = x2$ &   large OCV    &  $0 \le x1 \le 2$  &   12   & $55 k\text{\texteuro} $ & 0.8 & $40 \text{\ }t/day$\\
$F3 = x3$ &   barge    &   $0 \le x1 \le 2$   &   16 & $35 k\text{\texteuro}$  & 0.5 & $35 \text{\ }t/day$\\
\bottomrule
\end{tabular}}
\label{tab:owf_variables}
\end{table}

\noindent \textit{Project activities}\\
\noindent The installation process involves strategically placing and securing mooring anchors using a combination of vessels, adjusted for efficiency based on available deck space and potential risk event delays. Two distinct activities can be defined: (1) installing the available anchors on deck of the vessels; and (2) bunkering, which occurs after the installation sequence once all available anchors have been installed. Bunkering refers to the transfer of additional anchors to the vessels when the condition that the number of anchors left to install is zero or less is met.

The duration of these two activities are defined in \autoref{tab:c1_activites} and represented by the pessimistic $a_i$, most-likely $m_i$, and  optimistic $b_i$ time estimates. These values are usually obtained either from experience or past data of similar projects and make up the three-point estimates of the activities duration's, which are used to define the Beta-PERT distributions (see \autoref{methodology}).

\begin{table}[H]
\caption{Strategic planning application: Project activities}
{\begin{tabular}{l c c c}
\toprule
\multicolumn{1}{l}{Activity description} &
\multicolumn{3}{c}{Activity Duration (days)} \\
\cmidrule{2-4}
& \multicolumn{1}{l}{Optimistic $a_i$} &
\multicolumn{1}{l}{Most-likely $m_i$}  &
\multicolumn{1}{l}{Pessimistic $c_i$} \\
\midrule
Installation&   0.80    &   1.00   &   1.20    \\
Bunkering small OCV &   1.20    &  1.50  &   1.80    \\
Bunkering large OCV &   1.60    &   2.00   &   1.40    \\
Bunkering barge &   2.00    &   2.50   &   3.00    \\
\bottomrule
\end{tabular}}
\label{tab:c1_activites}
\end{table}

\noindent \textit{Risk events}\\
\noindent The offshore installation sequence increases in complexity due to uncertainties in the execution. External weather impacting the workability window of vessels, and project specific constraint often disrupt the execution, making completion according to initial forecast difficult to achieve. In probabilistic planning, risk events are included by considering their occurrence probability. \autoref{tab:risk_events} defines the risk events identified with their corresponding three-point duration estimates. The probability density $f(z_i; a_i, m_i, b_i)$ can be built using the three-point estimate. For every risk event, an occurrence probability $p_r$ is defined.

\begin{table}[H]
\caption{Strategic planning application: Risk events.}
\resizebox{\textwidth}{!}{\begin{tabular}{l c c c c c}
\toprule
\multicolumn{1}{l}{Risk Event} &
\multicolumn{3}{c}{Risk duration [$days$]} &
\multicolumn{1}{c}{Affected Activity} &
\multicolumn{1}{c}{Probability $p_r$} \\
\cmidrule{2-4}
& \multicolumn{1}{c}{$a_i$} &
\multicolumn{1}{c}{$m_i$}  &
\multicolumn{1}{c}{$b_i$} & & \\
\midrule
Weather Delay & 0.50 & 1.00 & 1.50 & Installation & 0.20 \\
Supply Chain Issue & 1.00 & 1.50 & 2.00 & Bunkering & 0.10 \\
Technical Failure & 0.50 & 1.00 & 1.50 & Installation & 0.15 \\
Logistical Constraints & 1.00 & 1.50 & 3.00 & Bunkering & 0.20 \\
Marine Traffic & 0.50 & 0.75 & 2.00 & Installation & 0.05 \\
Environmental Restrictions & 2.00 & 3.50 & 5.00 & Bunkering & 0.10 \\
Lack of Skilled Personnel & 1.00 & 2.00 & 3.00 & Installation & 0.08 \\
Equipment Shortage & 1.00 & 2.00 & 4.00 & Installation & 0.12 \\
\bottomrule
\end{tabular}}
\label{tab:risk_events}
\end{table}

\noindent The installation sequence (network) is represented and established using a discrete event simulation (DES), which depends on the planning variables, activities and risks defined above. The probabilistic network with the underlying logical links is expressed as follows:

\begin{equation}
    \mathcal{N}(x_1, x_2, x_3, \mathbf{y})
\end{equation}

\noindent where the variables $\mathbf{y}$ contain the information on the respective vessel deck space or are related to the risks (and duration of the activity).
To ensure that the project is executed by at least one vessel, the sum of all vessels on the project must be greater than one. The technical feasibility is guaranteed by a constraint ensuring the anchor’s resistance meets or exceeds the applied forces defined by further performance functions ($F_4$ ,..., $F_7$). For a complete exposition of these formulations and their underlying principles, the reader is referred to \cite{VanHeukelum2024}. 

\subsubsection{Objective functions}

\noindent The optimisation framework considers the following four objectives that form the link between the network performance function and the preference functions. For a complete exposition of these formulations and their underlying principles, the reader is referred to \cite{Wolfert2023} or \cite{VanHeukelum2024}. The following objective functions are rewritten based upon to suit the needs of the Odycon framework.\\

\noindent \textit{Objective 1: Project duration}\\
The project duration $\Delta$ extracted from the network performance function expressed as:
    
    \begin{equation}
        O_1 = O_{PD} = \Delta(\mathcal{N}(x_1, x_2, x_3, \mathbf{y}))
    \end{equation}

\noindent \textit{Objective 2: Installation cost}\\
The project cost depends on the day rates of the vessels $R_i$ multiplied with their active duration on the project $t_i$ given the network performance function, and the cost of the installed anchors, expressed as:
    
    \begin{equation}
        O_2 = O_{C}=(815M_a+40,000)n_a+\sum_{i=1}^3 x_i R_i t_i(\mathcal{N}(x_1, x_2, x_3, \mathbf{y}))
    \end{equation}
    
\noindent \textit{Objective 3: Fleet utilisation}\\
Fleet utilisation is represented by the probability of a vessel being better utilised in another project and is defined as:
    
    \begin{equation}
        O_3 = O_{F}= \prod\limits_{i=1}^3 p_i^{x_i}
    \end{equation}

\noindent \textit{Objective 4: Sustainability (CO\textsubscript{2} emissions)}\\
Emissions are defined as the daily vessel emissions $E_i$ multiplied with their active duration on the project $t_i$ given the network performance function:
    
    \begin{equation}
        O_4 = O_{S}=\sum_{i=1}^3 x_iE_it_i(\mathcal{N}(x_1, x_2, x_3, \mathbf{y}))
    \end{equation}

\subsubsection{Preference functions}
\noindent The preference functions for this demonstrative application are defined based on the input from industry project management professionals. The maximum preference ($P_{1,PD} = 100$) towards the objective $O_{PD}$ is defined as the target duration $T_{tar} = 90$ days. A Beta-PERT distribution is employed to model the preference towards project duration. This unique preference modelling reflects that, in principle, the project manager in reality is somewhat interested in an earlier project delivery than the exact target duration, but certainly not much later than this duration with a bit of slack. If there is a greater interest to deliver earlier, then the area under the preference function should be increased so that a more asymmetric function is resulting. This approach enables a continuous evaluation of preferences towards the project duration, facilitating more informed and balanced management decisions amidst uncertainty. The other three preference functions for cost, fleet utilisation and CO\textsubscript{2} emissions are respectively convex or concave functions with the max preference value ($P_{1,C} =  P_{2,F} = P_{2,S} =100$) for the lowest respective objective values (\textit{min} $O_i$) and with the lowest preference value ($P_{1,C} =  P_{2,F} = P_{2,S} =0$) for the highest possible objective values (\textit{max} $O_i$) respectively. The resulting preference functions, describing the relations between the different preferences and the objectives, are depicted in \autoref{fig:pref_ofw}. \\

\subsubsection{Results}

\noindent To retrieve outcomes, we will first have to estimate the weights to generate IMAP solutions as part of the simulation \& optimisation framework. The global weights $w_k$ are split equally between the energy provider ($w_1 = 0.5$) and the marine contractor ($w_2 = 0.5$). To reflect the objectives of their interest, their individual local weights are respectively: $w_{1,PD} = 0.60$ for project duration, $w_{1,S} = 0.40$ for CO\textsubscript{2} emissions, $w_{2,C} = 0.70$ for the installation costs, and $w_{2,F} = 0.30$ for fleet utilisation. \autoref{tab:ofw_weights} defines the weights of each of the preference functions.

\begin{table}[H]
\centering
\caption{Weights for each of the preference functions, according to $w'_{k,i} = w_k \cdot w_{k,i}$. Note, stakeholders can also reflect interest to all four objectives while the global weight distribution is constant.}
    {\begin{tabular}{lccccc}
    \toprule
    Stakeholder $k$ & $w'_{k,PD}$ & $w'_{k,S}$ & $w'_{k,C}$ & $w'_{k,F}$ & $w_k$ \\
    \midrule
    Energy provider   & 0.30 & 0.20 & 0.00& 0.00& 0.50\\
    Marine contractor  & 0.00 & 0.00& 0.35 & 0.15 & 0.50\\
    \midrule
    Total & 0.30 & 0.20 & 0.35 & 0.15 & \textbf{1.00} \\
    \bottomrule
    \end{tabular}}
    
    \label{tab:ofw_weights}
\end{table}

\begin{figure}[H]
    \centering
    \includegraphics[width=\linewidth]{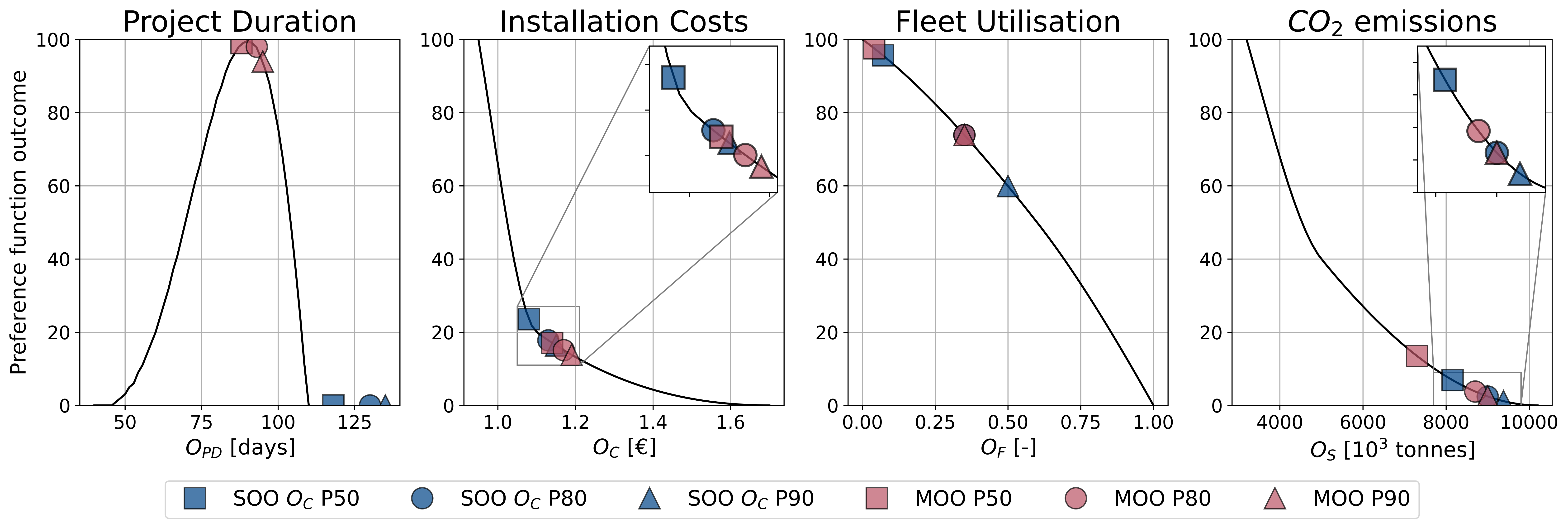}
    \caption{The four preference functions ($P_{1,PD}, P_{1,C}, P_{2,F}, P_{2,S}$) with respect to objectives ($O_{PD}, O_{C}, O_F, O_S$), including the preference score of the different optimisations and percentiles. The numerical results can be found in \autoref{tab:eval_ofw}.}
    \label{fig:pref_ofw}
\end{figure}

\begin{figure}[H]
    \centering
    \includegraphics[width=0.99\linewidth]{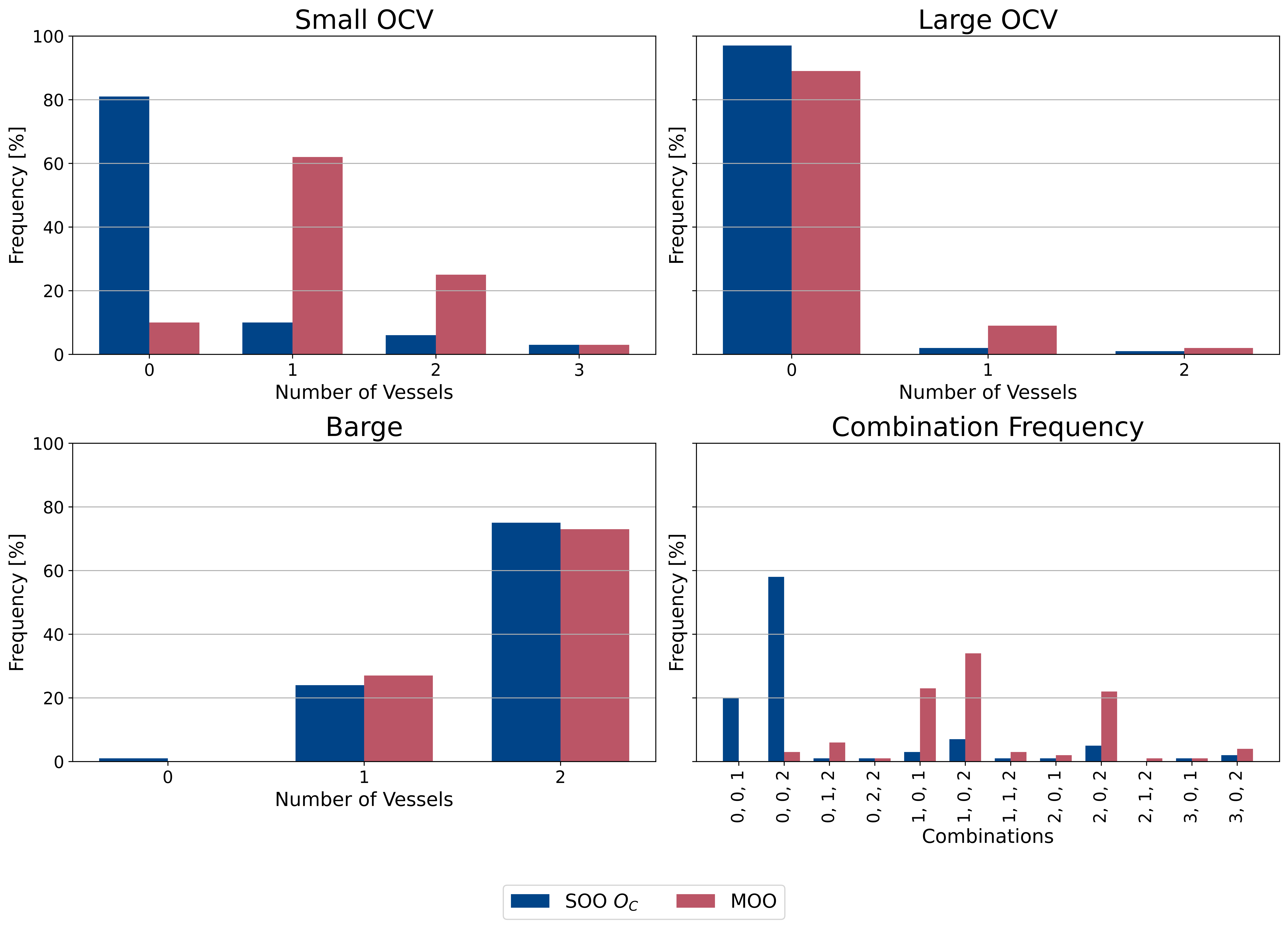}
    \caption{Criticality index of number of small OCV's (top left), large OCV's (top right), barges (bottom left) and vessel combinations (bottom right) for both optimisation cases. Note, only combinations that stored during the simulation \& optimisation are considered in this figure.}
    \label{fig:vessel_frequency}
\end{figure}

The final outcomes of the Odycon framework, are displayed in the form of frequency distribution diagrams referred to as criticality index. Within every MCS iteration, a new set of vessel combinations (planning variables) is derived through the IMAP optimisation. Given the unique realisation of the network, the vessel combinations change from one iteration to another. To identify the most critical set of planning variables, \autoref{fig:vessel_frequency} gives insight into the frequency of number of vessels per vessel type and the frequency distribution of optimal planning variables during the optimisation for both evaluation cases. The outcome reflects the percentage of MCS iterations of which a specific set of vessel combinations is derived. \autoref{tab:eval_ofw} compares the numerical results for the different evaluations given defined target-percentiles. The outcomes of the different percentiles per optimisation method are plotted in the different preference functions showing the respective preference given the percentile (see \autoref{fig:pref_ofw}). Note, as the MCS is stochastic in nature, it is not possible to provide a single (deterministic) outcome.\\

\begin{table}[H]
\caption{Evaluation of cases with 2000 MCS iterations ($O_{PD}$ in days, $O_C$ in $10^7$ \texteuro, $O_F$ in number of vessels [-], $O_{S}$ in tonnes). Where the MOO case considers $w'_{1,PD} = 0.30, w'_{1,S} = 0.20, w'_{2,C} = 0.35, w'_{2,F} = 0.15$.}
\resizebox{\textwidth}{!}{\begin{tabular}{l c c c c c c c c c c c c c}
\toprule
\multicolumn{1}{l}{Case} & \multicolumn{4}{c}{P50} & \multicolumn{4}{c}{P80} & \multicolumn{4}{c}{P90} \\
\cmidrule(lr){2-5} \cmidrule(lr){6-9} \cmidrule(lr){10-13}
& $O_{PD}$ & $O_C$ & $O_F$ & $O_{S}$ & $O_{PD}$ & $O_C$ & $O_F$ & $O_{S}$ & $O_{PD}$ & $O_C$ & $O_F$ & $O_{S}$ \\
\midrule
SOO $O_C$& 118 & 1.08 & 0.07 & 8152 & 130 & 1.13 & 0.35 & 9000 & 135 & 1.15 & 0.50 & 9380 \\
MOO & 88 & 1.14 & 0.04 & 7300 & 93 & 1.17 & 0.35 & 8700 & 95 & 1.19 & 0.35 & 9000 \\
\midrule
abs. diff.  & 30 & -0.06 & 0.03 & 852 & 37 & -0.04 & 0.00 & 300 & 40 & -0.04 & 0.15 & 380 \\
\midrule
rel. diff. & 25.4 & -5.6 & 42.9 & 10.5 & 28.5 & -3.5 & 0.0 & 3.3 & 29.6 & -3.5 & 30.0 & 4.0 \\
\bottomrule
\end{tabular}}
\label{tab:eval_ofw}
\end{table}

\noindent For comparative evaluation purposes, the MOO with the earlier defined weights and preferences is compared to a traditional SOO approach towards minimal cost (SOO $O_C$). The following conclusions can be made from the outcome:

\begin{enumerate}
    \item The SOO $O_C$ approach results in slightly lower costs for all percentiles (avg. 4\% lower) compared to the MOO approach by utilising no small OCV in 80 \% of the simulations, suggesting the cost-efficiency of the small OCV to be low (see \autoref{tab:eval_ofw}). It becomes evident that more vessels lead to higher costs.
    \item Comparing the SOO $O_C$ with the MOO approach, it becomes evident that considering project cost alone will not result in configurations that lead to overall favourable outcomes. While the project costs differ very little between the approaches, the MOO outperforms the SOO $O_C$ for all other objectives and percentiles leading to lower project duration's (avg. 28\% lower), lower probability of the fleet being better utilised elsewhere (avg. 24\% lower), and lower emission rates (avg. 6\% lower). This can be achieved by introducing small OCVs. This is explained by the reflection of interests towards the sustainability and fleet utilisation in the MOO. This comparison highlights that focusing solely on cost (single objective approach) does not accurately capture the complexities of stakeholder-oriented project planning. In contrast, utilising the MOO approach offers a more balanced perspective by taking into account all stakeholder interests and achieves a higher aggregated preference compare to the SOO $O_C$ approach, including the inherent uncertainties and stochastic nature of the project leading to overall favourable project realisations.
    \item \autoref{fig:pref_ofw} illustrates that as the percentile increases, the preference for both cases decreases concerning their respective objectives. This trend indicates that decision-makers encounter lower preference outcomes when aiming for higher objective probabilities. The target probability can thus not only be reflected by the unit of the individual objectives but also directly in the common preference domain. By examining the target-preference values, decision-makers can identify the probability of achieving a particular outcome.
    \item The criticality index (see \autoref{fig:vessel_frequency}) aids decision-making by clarifying the utilisation of specific types of vessels and their combinations, thereby identifying the most effective vessel configurations given the underlying uncertainties and stochastic factors. For instance, in both optimisation scenarios, large OCVs are not utilised in 80\% of the simulations. These observations provide project managers with actionable insights, enabling them to allocate their planning and tendering budgets more effectively towards solutions that are both efficient and feasible.
\end{enumerate}

\subsection{Dynamic Control application: Adapt Flexibly}

\noindent\textit{Technical context:}
The Schiphol-Amsterdam-Almere (SAA) Project is the largest road construction project in the Netherlands in the period from 2012 to 2024. It realises an expansion of the capacity of the Dutch national highways A6, A1, A10-east and A9 and thus aims to improve the accessibility of the northern part of the Randstad (Schiphol, Amsterdam, Almere). Project SAA includes 63 kilometres of motorway, the construction of two tunnels, two large bridges and an aqueduct; the modification of five interchanges, and about 100 engineering structures. As part of this project, the SAAone consortium (Boskalis, Volker-Wessels, Hochtief and DIF), together with the Dutch highways agency Rijkswaterstaat, provided the design, construction, maintenance and operation. The major construction project took more than four years and consisted of a massive schedule that could be regularly adjusted on the run. \\

\noindent\textit{Social context:}
This application includes the following concurrent stakeholder objectives both from SAAone and RWS: (1) meeting target completion, (2) managing project control costs, and (3) minimising the traffic nuisance. In order to manage the project in the best possible way, the costs of the extra control measures to be allocated, which are possibly not part of SAAone's risk budget, are at odds with the opening of the new highway on time (interest RWS). Even more, steering solely on the completion date from RWS's project delivery department perspective may also be of opposite interest to minimising traffic nuisance during construction, which RWS's service operations department is responsible for. That is what makes this example a dynamic multi-objective control and tripartite stakeholder problem.

We will now first describe the dynamic control optimisation problem by working through the mathematical statement (see \autoref{sec:optimisation}), resulting in performance-, objective-, and preference functions.

\subsubsection{Performance functions}

\noindent \textit{Project activities}\\
The project consists of a number of distinct activities, as detailed in \autoref{tab:case3 activity}. For each activity, the durations are presented by an optimal $a_i$, most likely $m_i$, and pessimistic $b_i$ time estimate. These values contribute to the three-point estimates for the activity duration's, forming the basis for the Beta-PERT distributions. Additionally, the sequence of activities is defined by specifying predecessors for each activity, thereby integrating the dependencies among activities into the network model.

Next to the individual uncertainty of activity duration's a set of shared uncertainties is considered to model activity duration correlations arising from common factors (e.g., weather, labour skills). The implementation of stochastic activity correlation was implemented in accordance with the methodology and underlining example presented in \cite{Kammouh2022}. For further information regarding the evaluation of networks under correlated uncertainty, see \cite{WangDemsetz2000}. \autoref{tab:case3 correlation} defines the uncertainty factors, their three-point estimates and relationships with the project activities.\\

\noindent \textit{Risk events}\\
The project is affected by many possible risk events that ad an additional source of uncertainty that can negatively affect the desired project outcome (see \autoref{tab:saa_risks}). The impact of these risks is defined by their risk duration (three-point duration estimates) in the same manner as done for the activities and control measures. Every risk is characterised by a given probability of occurrence ($p_e$) and the relations between the risk events and activities.\\

\noindent \textit{Control measures}\\
\noindent Based on the project database, a number of possible control measures $x_n$ are identified (see \autoref{tab:saa_measures}), that can increase the probability of finishing the project at the target completion time. These control measures, identified by the project manager, take into account limitations on material and human resources. These control measures are characterised by their impact towards the three identified objectives: (1) the number of days reducing the duration of the affected activities, (2) the cost of a measure, and (3) the impact on traffic nuisance, each given by a three-point estimate (minimum, most-likely, and maximum). The correlation factor ($\eta$) is defined for accurately determining mitigation effects in construction delays by accounting for the nonlinear relationship between mitigated duration and effects, varying from direct proportional effects for additional personnel or equipment to non-proportional one-off expenses.\\

\noindent The SAAone construction sequence (see \autoref{fig:gantt}), is established using a logical network model and build upon the planning project activities and risks, depending on the control measures. The probabilistic network with the underlying logical links is expressed as follows:

\begin{equation}
    \mathcal{N}(a_n *  c_{i,n}, \mathbf{y})
    \label{network_performace_fucntion_saa}
\end{equation}

As noted in remark 4 of \autoref{sec:optimisation} the methodology defines the implementation of a control measure  to be represented by allocation of a measure $a_n$ multiplied with its impact $c_{i,n}$ towards objective $i$ as following: $F_n = x_n = a_n * c_{i,n} $ (see the set of control measures in \autoref{tab:saa_measures}). \\

\subsubsection{Objective functions}
\noindent The optimisation framework considers the following four objectives that form the
link between the network performance function and the preference functions.\\

\noindent \textit{Objective 1: project duration}\\
The project duration $\Delta$ is extracted from the network performance function and can be expressed as follows:

\begin{equation}
    \label{eq:objective_time}
    O_1 = O_{PD} = \Delta(\mathcal{N}(a_n *  c_{1,n}, \mathbf{y}))
\end{equation}

\noindent where $O_{PD}$ is expressed in days.\\

\noindent \textit{Objective 2: project control cost}\\
The project control cost is defined as the sum of the cost impact of active control measures. As the demonstrative case is defined according to a \textit{Design, Build, Finance and Maintain-contract} (DBFM) contract, the objective function is extended by a penalty and reward scheme (contractual project completion performance scheme), it reads as:

\begin{equation}
    \label{eq:objective_cost}
    O_2 = O_C = \sum_{n=1}^{N} \left(a_n \cdot c_{2,n} \right) + \Delta_1 \cdot P_c - \Delta_2 \cdot R_c
\end{equation}

Subject to:

\begin{equation}
\label{constraint_1}
    \Delta_1 = 
\begin{cases} 
\Delta - T_{\text{tar}}, & \text{if } \Delta > T_{\text{tar}} \\
0, & \text{otherwise}
\end{cases}
\end{equation}

\begin{equation}
\label{constraint_2}
    \Delta_2 = 
\begin{cases} 
T_{\text{tar}} - \Delta, & \text{if } \Delta < T_{\text{tar}} \\
0, & \text{otherwise}
\end{cases}
\end{equation}

\noindent where $O_C$ is expressed in \texteuro, $\Delta_1$ is the project delay after implementing the control measures, $\Delta_2$ is the duration reduction beyond target duration, $ T_{tar}$ is the target (i.e., desired/planned) project duration specified by the involved stakeholders, $P_c$ is the daily penalty, and $R_c$ is the daily reward. It holds that a penalty and reward cannot occur simultaneously, thus $\Delta_1 * \Delta_2 = 0$. Note, the project duration after implementation of control measures $\Delta$ is extracted from the network performance function $\mathcal{N}(a_n *  c_{2,n}, \mathbf{y})$, see \autoref{network_performace_fucntion_saa}.\\

\noindent \textit{Objective 3: traffic nuisance}\\
During construction, the experience of a road user (considered as traffic nuisance) is heavily impacted. The impact on road user experience during construction is becoming an increasingly important aspect within infrastructure project execution. Currently, most service operations departments measures the effect on car traffic in terms of Lost Vehicle Hours (LVHs), which is then converted into a monetary value and considered in a contractual performance scheme. However, this approach overlooks the actual impact on the user which can not be monetised, respectively maximising their experience during construction becomes essential for improving service operations during construction. Similar to the penalty and reward scheme of the objective cost, a scheme is implemented to account for nuisance increase with the factor $P_e$ in case of a delay and decrease in traffic nuisance with the factor $I_e$ in case of early completion. $\Delta_1$ and $\Delta_2$ are considered according to \autoref{constraint_1} and \autoref{constraint_2}. Again, it holds that a penalty and reward cannot occur simultaneously, thus $\Delta_1 * \Delta_2 = 0$. The objective traffic nuisance reads as follows:

\begin{equation}
\label{o:experience}
O_3 = O_{N} = \left( \frac{1}{\sum_{n=1}^{N} \left(c_{3,n} \right)} \sum_{n=1}^{N} \left(a_n \cdot c_{3,n} \right) \right) * S + \Delta_1 * P_n - \Delta_2 * R_n
\end{equation}

where $S$ is the scaling factor of 10. Since there is no standard method on how to represent the traffic nuisance, a scale from 0 (baseline nuisance during construction) to 10 (worst possible traffic nuisance during construction) is introduced.

\subsubsection{Preference functions}
\noindent The preference functions for this project control application were developed with one of the co-authors who was involved in the project during the execution phase. In accordance with the project RWS delivery department and the SAAone consortium, the preference curve for project duration (Time) is defined by three duration estimates using a Beta-PERT distribution. The maximum preference ($P_{1,PD} = P_{2,PD} = 100$) towards objective $O_{PD}$ is set at the contract target duration, $T_{tar} = 1466$ days. The minimum preference ($P_{1,PD} = P_{2,PD} = 0$) is determined by two points: the target duration plus the maximum allowable delay, and a minimum duration of 966 days. Note, the preference towards the project delivery time is expressed by the joint contract of RWS and SAAone, therefore the one preference curve is defined towards the objective target duration. Similar to the prior example, the preference curve is established using a Beta-PERT distribution to accurately reflect the project manager's preferences for the project duration. This unique preference modelling reflects that, in principle, the project manager in reality is somewhat interested in an earlier project delivery than the exact target duration (blue vertical line), but certainly not much later than this duration with a bit of slack (red vertical line). If the project manager has a 'fixed' project delivery with no slack, the red and blue line should coincide, essentially reflecting no preference for any delay. If there is an increased interest to deliver earlier, then either the target duration can be changed or the area under preference function should be increased to create an extremely asymmetric function. Since the delay of the project changes with every iteration of the MCS, the preference curve is adapted accordingly for every realisation of the project network. The preferences for lowest project control cost and lowest nuisance are considered to have the highest value ($P=100$). Conversely, the lowest preference value ($P=0$) is considered for the total sum of all control measures' cost and for the worst possible level of traffic nuisance. For an initial estimate, both preference functions (cost and nuisance) are linear between these lowest and highest preference values. The three resulting preference functions, which describe the relations between different values for $P_{1..3,1..3}$ with respect to objectives $O_{1..3}$, are shown in \autoref{fig:pref_mico}.

\subsubsection{Results}

\noindent To retrieve Odycon outcomes, we will first have to estimate the weights to generate IMAP solutions as part of the MCS results. Within a concurrent discussion between the two (global) stakeholders RWS ($w_1 = 0.5$) and the SAAone consortium ($w_2 = 0.5$) the three individual weights are discussed to reflect their shared objectives. 
RWS internal interests are represented by the project delivery department with $w_{1,PD} = 0.70$ for the project duration and the service operations departments with $w_{1,N} = 0.30$ for traffic nuisance. SAAone expressed their individual (local) weights to be $w_{2,PD} = 0.70$ for project duration and $w_{2,C} = 0.30$ for the control costs. \autoref{tab:SAA_weights} reflects the resulting weights of each of the preference functions. As mentioned before, the preference towards the project delivery time is expressed by the joint contract of RWS and SAAone, therefore the one preference curve is defined towards the project duration.

\begin{table}
\centering
\caption{Weights for each of the preference functions, according to $w'_{k,i} = w_k \cdot w_{k,i}$. Note, stakeholders can also reflect interest to all objectives while the global weight distribution is constant.}
    {\begin{tabular}{lccccc}
    \toprule
    Stakeholder $k$ & $w'_{k,PD}$ & $w'_{k,C}$ & $w'_{k,N}$ & $w_k$ \\
    \midrule
    RWS & 0.35 & 0.00 & 0.15 & 0.50\\
    SAAone consortium & 0.35 & 0.15 & 0.00 & 0.50\\
    \midrule
    Total & 0.70 & 0.15 & 0.15  & \textbf{1.00} \\
    \bottomrule
    \end{tabular}}
    
    \label{tab:SAA_weights}
\end{table}

For comparative evaluation purposes, three scenarios are considered. First, there is a SOO focused on minimising costs (SOO $O_C$). Second, another SOO targets minimising nuisance (SOO $O_N$). Finally, a MOO is conducted using a defined set of weights (see \autoref{tab:SAA_weights}) to reflect all stakeholder's interests. To ensure convergence towards the project target date $T_{tar}$ a penalty of $P_c = 10k \text{\texteuro}$ and $P_n = 0.1$ is defined for SOO costs and SOO nuisance case respectively. Without these "constraints" the SOO would result in no used control measures due to the missing incentive to utilise them. This is however not needed for the MOO, as it solved the aforementioned need for direct time-cost trade-off techniques.

\begin{figure}[H]
    \centering
    \includegraphics[width=0.99\linewidth]{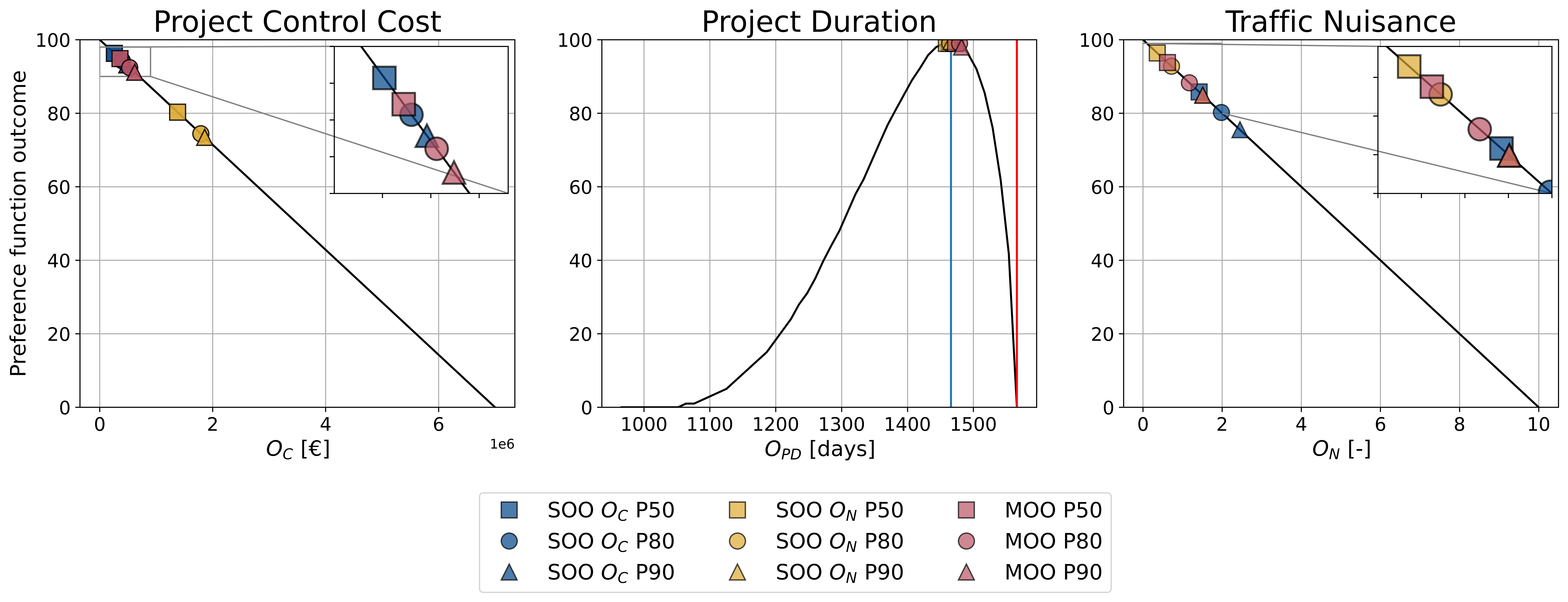}
    \caption{Preference functions ($P_{1, PD}, P_{2, C}, P_{1, N}$) with respect to objectives ($O_{PD}, O_{C}, O_{N}$), including the  preference score of the different optimisation scenarios and percentiles. The numerical results can be found in \autoref{tab:eval_saa}. The blue line represents the target duration $T_{tar} = 1466$ days and the red line an exemplary maximum delay of 100 days. Note that, as previously mentioned, since the project delay changes with each iteration of the MCS, the preference curve is adapted accordingly for each realisation of the project network. However, this adaptation cannot be visualised. }
    \label{fig:pref_mico}
\end{figure}

\begin{figure}[H]
    \centering
    \includegraphics[width=\linewidth]{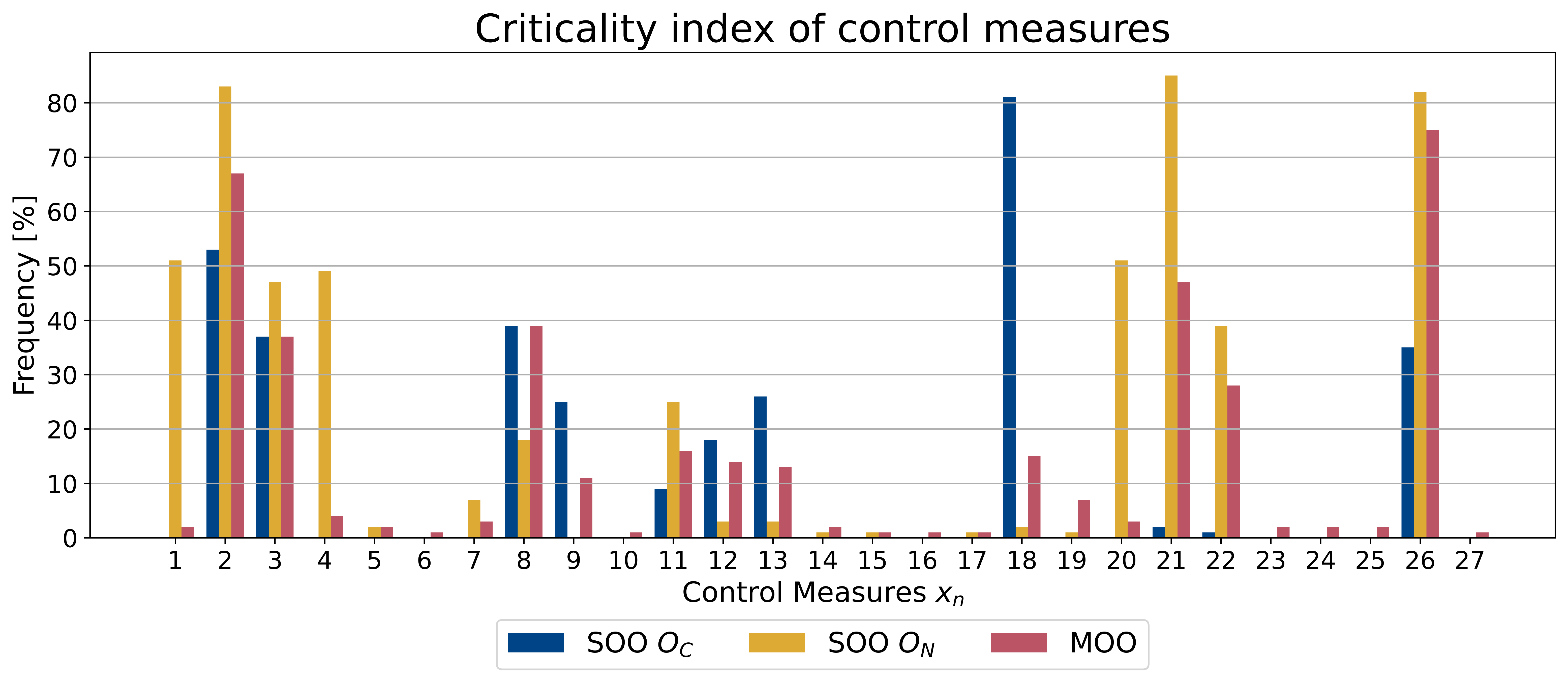}
    \caption{Criticality index of mitigation measures for all three optimisation scenarios.}
    \label{fig:saa_frequency}
\end{figure}

\begin{figure}[H]
    \centering
    \includegraphics[width=\linewidth]{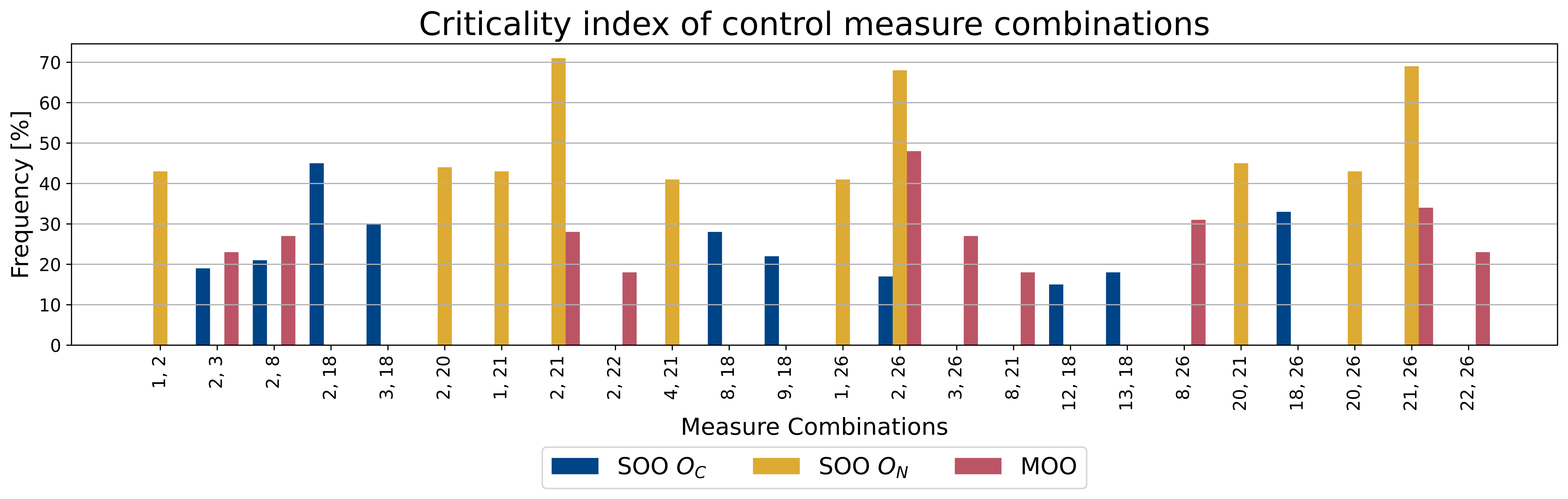}
    \caption{Criticality index of mitigation measures combinations for all three optimisation scenarios.}
    \label{fig:saa_frequency_combinations}
\end{figure}

The final outcomes of the Odycon framework (integration of IMAP and MCS) are displayed in the form of frequency distribution diagrams referred to as criticality index. Similar to the strategic planning application, the criticality index of control measures is derived. Within every MCS iteration a new set of control measures is derived. To identify the most critical set of control variables \autoref{fig:saa_frequency} gives insight into the frequency of control measures and \autoref{fig:saa_frequency_combinations} on the frequency of combined occurrence of control measures. \autoref{tab:eval_saa} compares the numerical results for the different evaluation cases given defined target-percentiles. The results of the different percentiles and cases are respectively displayed in the different preference functions showing the respective preference given the percentile (see \autoref{fig:pref_mico}).\\

\begin{table}[H]
\caption{Evaluation of cases with 2000 MCS iterations ($O_{PD}$ in days, $O_C$ in $10^6$ \texteuro, $O_N$ in nuisance level [0,10]). Where the MOO case considers $w'_{1,PD} = 0.70, w'_{2,C} = 0.15, w'_{1,N} = 0.15$.}
\resizebox{\textwidth}{!}{\begin{tabular}{l c c c c c c c c c c c c}
\toprule
\multicolumn{1}{l}{Case} & \multicolumn{3}{c}{P50} & \multicolumn{3}{c}{P80} & \multicolumn{3}{c}{P90} \\
\cmidrule(lr){2-4} \cmidrule(lr){5-7} \cmidrule(lr){8-10}
& $O_{PD}$ & $O_C$ & $O_N$ & $O_{PD}$ & $O_C$ & $O_N$ & $O_{PD}$ & $O_C$ & $O_N$ \\
\midrule
SOO $O_C$ (C) & 1463 & 0.26 & 1.42 & 1468 & 0.40 & 1.98 & 1471 & 0.48 & 2.45 \\
SOO $O_N$ (N) & 1459 & 1.38 & 0.36 & 1464 & 1.79 & 0.72 & 1466 & 1.86 & 0.98 \\
MOO (M)& 1473 & 0.36 & 0.62 & 1479 & 0.53 & 1.17 & 1482 & 0.62 & 1.51 \\
\midrule
abs. diff. (C - N) & 4 & -1.12 & 1.06 & 4 & -1.39 & 1.26 & 5 & -1.38 & 1.47 \\
abs. diff. (C - M) & -10 & -0.10 & 0.80 & -11 & -0.13 & 0.81 & -11 & -0.14 & 0.94 \\
abs. diff. (N - M) & -14 & 1.02 & -0.26 & -15 & 1.26 & -0.45 & -16 & 1.24 & -0.53 \\
\midrule
rel. diff. (C - N) & 0.3 & -430.8 & 74.6 & 0.3 & -347.5 & 63.6 & 0.3 & -287.5 & 60.0 \\
rel. diff. (C - M) & -0.7 & -38.5 & 56.3 & -0.7 & -32.5 & 40.9 & -0.7 & -29.2 & 38.4 \\
rel. diff. (N - M) & -1.0 & 73.9 & -72.2 & -1.0 & 70.4 & -62.5 & -1.1 & 66.7 & -54.1 \\
\bottomrule
\end{tabular}}
\label{tab:eval_saa}
\end{table}

\noindent The following conclusions can be made from the comparative evaluation outcome:

\begin{enumerate}
    \item Comparing the three cases it becomes evident that for all percentiles the relative difference in the project duration ($O_{PD}$) is very little (1\% difference between cases) suggesting that all cases perform well with optimising towards meeting the target duration. However, the SOO $O_C$ and SOO $O_N$ are only able to achieve this by the implemented contractual penalty (considered as a trade-off technique). The outcome is thus highly influenced by the defined penalty, limiting the reflection of actual project reality but instead only the contractual "pseudo-reality'. The MOO approach achieves similar results given the desired preference of all involved stakeholders allowing for an independent representation towards the project duration and control cost, overcoming the need for monetisation techniques. Moreover, these monetisation techniques are limited with regard to "soft" objectives (as traffic nuisance) where direct cost relation and quantification is not possible.

    \item As expected, larger relative differences in case of project planning costs and traffic nuisance are visible. The SOO $O_C$ results in an average of 33\% lower cost compared to the MOO approach. The SOO $O_N$ results in an average of 63\% lower traffic nuisance compared to the MOO approach. However, when considering a stakeholder-oriented behaviour, it becomes evident that the two SOO cases reflects a narrow and one-sided view not accounting for other set of project goals. With the SOO approaches no synthesis that benefits the overall project outcome is achieved, revealing a significant drawback. The MOO outcome however resulting in acceptable levels of all respective objectives and the overall highest aggregated preference for the respective percentiles. This clearly demonstrates that the synthesis towards all project objectives can result in overall group satisfaction. Similar to the strategic planning application, utilising the MOO approach offers a more balanced perspective by taking into account multiple objectives, including the inherent uncertainties and stochastic nature of the project leading to overall favourable project realisations. Similar to the prior example it becomes evident in \autoref{fig:pref_mico} that with increasing percentile, the preference decreases concerning their respective objectives. As expected, this trend indicates that decision-makers face lower preference outcomes when targeting higher objective probabilities.
    
    \item The criticality index (see \autoref{fig:saa_frequency} and \autoref{fig:saa_frequency_combinations}) aids decision-making by clarifying the utilisation of specific control measures and their combinations, thereby identifying the effectiveness of certain control measures and combinations given the underlying uncertainties and stochastic factors. Comparing the three cases, the criticality index differs significantly where certain measures seem to have clear positive impact towards certain objectives. For example measure 18 is utilised in 80\% of the iterations in the SOO $O_C$ however shows little impact for the SOO $O_N$ or MOO case. The results suggest that some corrective measures (i.e. with low occurrence frequency, e.g.: 14, 15, 16, 17) can be ignored, while others should be prioritised. The frequency of measures given SOO $O_N$ and SOO $O_N$ do not consider the impact a certain measure has towards other objectives, neglecting associative decision-making. The MOO approach reflects the frequency of measures that have the highest aggregated group preference towards all objectives revealing their actual impact. However, to arrive at effective project control strategies, it is important not to rely solely on the most critical corrective measure as a standalone solution. Instead, these measures should be combined to form a comprehensive mitigation strategy selected by the project manager. To this extent, the criticality index of measure combinations gives insight into the correlated impact of measures. The likelihood of timely project completion should then be verified through a simulation allocating these specific measures. Similar critically analysis can be performed for the project activities and project paths (see \cite{Kammouh2021}).

\end{enumerate}

\noindent These observations provide project managers with actionable insights, enabling them to define the most effective project control measures given all the project objectives. A dynamic adaption and verification thought the project stages is enabled with the ability to adapt the set of feasible measures.

\section{Discussion \& next steps for further development}
\noindent The two example application demonstrate the advances and added benefit of the introduced Odycon methodology. Facilitating associative decision-making in uncertain projects come with limitations towards the technical and social integration. The challenge lies in the combined complexity of synthesising the concurrent interests of all stakeholders while managing intricate technical interdependencies effectively.
The methodology is built for computer-aided and data-driven decision-making, relying heavily on the availability and accessibility of large project performance information and stakeholder preferences mapping. For example, estimates for activity durations, technical performance data of engineering assets, occurrence probabilities of risks, and data on external environmental influences significantly impact the interdependency and workability of a project. To achieve reproducible and consistent outputs, the variability in data quality from different sources must be minimal. Moreover, accurately representing stakeholder preferences and their relative individual weights significantly impacts Odycon's outcomes. Consequently, well-informed project decision-making is fundamentally limited by the availability and quality of both preference and performance information.

While the introduced Odycon methodology is adaptable to various project specifics and constraints, the inflexibility of existing contracts may hinder the accommodation of dynamic changes. Legal constraints and compliance issues could limit the effectiveness of this open and cooperative decision-support methodology, impeding adaptive and associative project management. 

Towards that extent, the following further research and development is recommended:

\begin{enumerate}

    \item As demonstrated in the example applications, the results depends on a good reflection of human objectives and their preferences \& weights (social integration). Further research is needed towards establishing a Structured Stakeholder Judgement model to improve and refine the preference elicitation to estimating the preference functions and weight distributions. This type of model could improve cooperative decision-making, promoting a collective mindset shift among all involved parties to fully utilize the benefits of the Odycon methods for the advantage of all project associated stakeholders. As a first step, the so-called choice-based conjoint analysis method can be introduced to generate an initial estimate for this.
    
    \item  To improve the accurate reflection of uncertainty (technical integration) the implementation of Structured Expert Judgment (e.g. Cooke’s model) should be considered when decision-makers are confronted with insufficient data availability \citep{Leontaris2019}. In case of large sets of available field-data, the use of real-time data methods like non-parametric Bayesian Networks could be integrated to model uncertainty without relying on predefined distributions \citep{leontaris2018probabilistic}.
    
    \item The example applications considered Beta-PERT modelling, neglecting other simulation-based models. To improve the simulation of the network, the impact of models that consider the integration of different types of activity links and automated changes in the network structure should be explored \citep{wang2005impact}.
    
    \item Constraints regarding resources limitation (e.g., personnel, material, assets) are currently neglected. However, for more realistic applications, the availability of certain resource pools during the project execution should be incorporated to give project managers a more accurate representation of feasible project planning and control strategies.

    \item Both postmortem analysis and real life projects research can further support the effect of using Odycon and add confidence to contractor-client relationships so that the current rigid forms of contract may become more open. 

    \item Further research can be done to use Odycon in stalemate projects as a transparent and objective mediator to 'confront project conflicts into yes'.

    \item The Odycon methodology can be further extended for project portfolio optimisation to reduce inefficiencies. 

\end{enumerate}

\section{Conclusion}
\noindent The increasing complexity of projects comes with the need for stakeholder and goal-oriented decision-making, and is addressed by following a systems-oriented approach. 

This paper presents Odycon, a pure \textit{a-priori} stochastic simulation \& optimisation methodology integrating the capability of the project (technical domain), the human goal-oriented behaviour (human domain), and the association of stakeholder-oriented behaviour (social domain). To this extent, the IMAP optimisation method is integrated into a probabilistic MCS, accurately reflecting project uncertainties while enabling best-fit for common-purpose project management. The advances towards strategic planning and dynamic control were demonstrated with two distinct examples demonstrating the utilisation during up-front planning and on-the-run project control. The following general conclusions can be drawn.

First, the complex reality of project management is addressed by identification of the most effective and desired set of project planning \& control variables given the uncertainties and stochastic factors while considering overall group satisfaction. The interplay with stakeholder-oriented concurrent objectives is addressed by the integration of the ’associative preference domain', enabling an optimal solution that best fits the common purpose, using \emph{a-priori} optimisation. This approach supports satisfactory project planning and control, allowing stakeholders to meet their individual needs through collective performance. The model results of both demonstrative applications show that joint project success, in alignment with the defined common purpose, will increase when individual stakeholders relinquish their pure self-interest, ultimately enabling an optimal solution that benefits the whole group. Here, Odycon aligns with previous empirical results of management studies, as mentioned in the introduction. The developed methodology allows for the integration and combination of variables from different project phases, as well as flexible adaptation and development tailored to project and stakeholder specific needs. 

Second, Odycon opens the possibility to reflect all types of project goals, with no need for cost or resource trade-off techniques, thus providing flexibility and better insight into different objectives and system behaviour. It enables adaption towards project specific circumstances where numerous important factors may not be covered by contractual specifications or budget, leading to a decline in quality under the conventional management practice. This consideration is essential for enhancing project planning and control strategies while maintaining project integrity in terms of quality, budget, and timeline. While the budget of a project is determined by its financial limits, other objectives describing the quality or impact of a project should be evaluated based on their degree of 'satisfaction,' which reflects the utility or value they offer.

Third, Odycon provides clear guidance in budget and resource allocation before and during project execution, making it an essential tool to support decision-making. Project managers and associated stakeholders are able to jointly create actionable insights into the system behaviour and decision-making impact of planning and/or control variables given the underlying uncertainties a project faces. This enables to jointly allocate planning and control budgets more effectively towards solutions that are both feasible and desired.

Fourth, Odycon enables well-informed and deliberative decision-making by removing bias from the decision-making process. Manual trial-and-error approaches are time-consuming and inefficient, relying on 'gut instinct' and human judgment on uncertainties, which can be biased and limited. In contrast, Odycon's computer-aided decision-support framework uses simulation and optimisation models to explore numerous solutions systematically, reducing human biases, enhancing efficiency, and acknowledging solution spaces no longer conceivable for humans. Odycon thus facilitates a shift within complex project management from unsubstantiated trial-and-error to transparent and data-driven decision modelling, optimizing over large numbers of variables and objectives, within a multidimensional solution space.\\

\noindent With this methodology, Odycon takes a next step in computer-aided data-driven decision-making for project and operations management through significantly improving the efficiency and effectiveness of decision-making, thereby enabling an associative ideal within reach.

\section*{Disclosure statement}
\noindent The authors report that there are no competing interests to declare.

\section*{Data availability statement}\label{data availability}
\noindent The two demonstrative applications for strategic project planning and dynamic project control are available via a GitHub repository online in accordance with confidentiality agreements. The repository is located here: \url{https://github.com/Boskalis-python/ODYCON}. Note, for preference aggregation as part of the optimisation, the following algorithm is used: \url{https://github.com/Boskalis-python/a-fine-aggregator}.

\section*{Acknowledgements}
The authors would like to thank Dr. Ruud Binnekamp of Delft University of Technology for his valuable reflections and guidance in completing this paper.

\appendix

\section{Preference score aggregation}\label{a-fine-a}

\noindent The following describes the aggregation algorithm for use in preference-based design and decision-making. In short, the following two starting principles apply to this algorithm (see \citet{Barzilai2022} or \citet{Wolfert2023}):

\begin{enumerate}
    \item it should reflect relative scoring as encountered in actual design and decision-making practice.

    \item it should adhere to the governing mathematics in a one-dimensional affine space, which is the mathematical model applicable to preference score(s).
\end{enumerate}

The algorithm therefore consists of two operations: (1) normalising the preference scores of all alternatives per criterion, and (2) finding the representative aggregated preference score $P^*$ for each alternative using the weighted least squares method. These two operations are further elaborated mathematically below.

\subsection{Normalisation}

\noindent For normalisation, the standard score (z-score) method is used. This yields a normalisation that preserves information about the population of preference scores and reads as follows:

\begin{equation}
    z_{i,j} = \frac{p_{i,j}-\mu_j}{\sigma_j}\textrm{ for }
    i=1,2,...,I \textrm{; }
    j=1,2,...,J
\end{equation}

Here $z_{i,j}$ is the normalised score of alternative $i$ for criterion $j$; $p_{i,j}$ is the preference score of alternative $i$ for criterion $j$; $\mu_j$ is the mean of all preference scores $p$ for criterion $j$; $\sigma_j$ is the standard deviation of all preference scores $p$ for criterion $j$. By performing this normalisation for all criteria $J$, the preference scores are transformed to a single scale with the same properties ($\mu_J=0$, $\sigma_J=1$). 

\subsection{Weighted least squares}

\noindent Since all $z_{i,j}$ scores are now on a single scale, it is possible to compare all normalised scores per alternative with each other. To find the representative aggregated preference score of an alternative that provides a best fit of all (weighted and relative) scores for each criterion, a minimisation of the weighted least squares difference between this aggregated score and each of the (normalised) individual scores on all criteria is applied. This is expressed mathematically as follows:

\begin{equation}
     \textit{Minimise} \text{\ }S_i = \sum_{j=1}^{J} w_j * (z_{i,j} - P_i^*)^2
\end{equation}

Note that since the search is for a single representative aggregated preference score, the model function $f(x_{i,j}, \beta_i)$ from the classical weighted least square method is replaced by $P_i^*$. The solution to this minimisation can be found by differentiating with respect to $P_i^*$ and equating it to zero. Since $\sum_{j=1}^{J}w_j=1$, this results in the following analytical expression for the representative aggregated preference score:

\begin{equation}
    P_i^* = \sum_{j=1}^{J} w_j * z_{i,j}
\end{equation}

\section{Dynamic control application data}\label{control data}

\begin{figure}[H]
    \centering
    \includegraphics[width=1.0\linewidth]{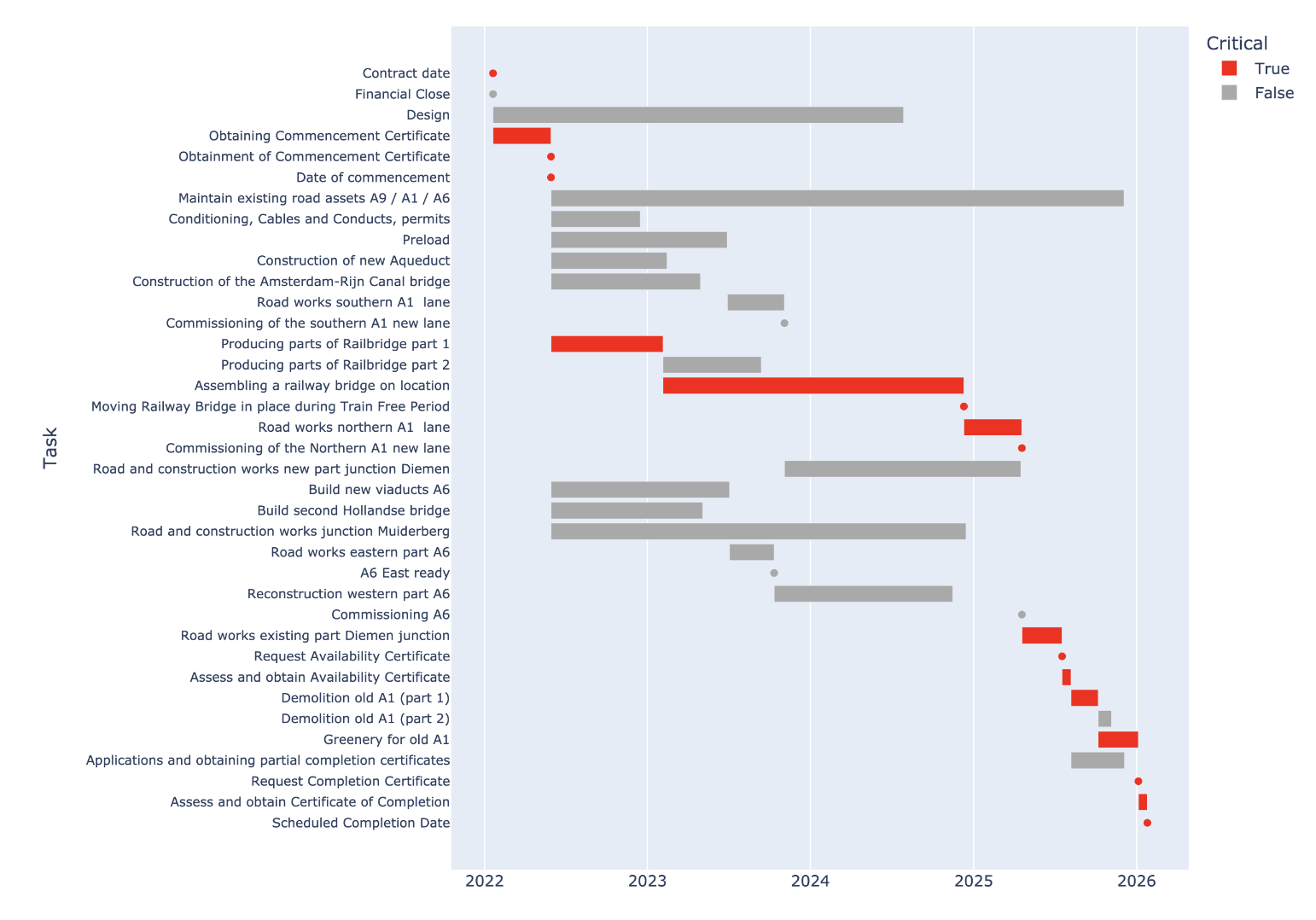}
    \caption{Simplified project planning and SAAone's construction sequence (logical network model) represented as a Gantt chart, a (deterministic) critical path is represented in red. This was also assumed as the base project planning in \cite{Kammouh2022}.} 
    \label{fig:gantt}
\end{figure}

\begin{table}[H]
    \centering
    \caption{Dynamic control application: Project Activities.}
    \label{tab:case3 activity}
    \resizebox{\textwidth}{!}{
        \begin{tabular}{clrrrr}
\hline \multirow[b]{2}{*}{ ID } & \multirow[b]{2}{*}{ Activity description } & \multicolumn{3}{c}{ Activity duration (days) } & \multirow[b]{2}{*}{ Predecessors } \\
\hline & & $a_i$& $m_i$&  $b_i$& \\
\hline
1 & Contract date & 0 & 0 & 0 & 0 \\
2 & Financial Close & 0 & 0 & 0 & 0 \\
3 & Design & 819 & 920 & 1435 & 1 \\
4 & Acquiring the certificate of commencement & 105 & 130 & 194 & 1 \\
5 & Commencement certificate is acquired & 0 &
0 & 0 & 2,4 \\
6 & Date of commencement & 0 & 0 & 0 & 5 \\
7 & Maintain existing road assets A9 / A1 / A6 & 976 & 1284 & 1836 & 6 \\
8 & Conditioning, Cables and Conducts, permits & 168 & 200 & 268 & 6 \\
9 & Preload & 324 & 395 & 525 & 6 \\
10 & Constructing a new Aqueduct & 200 & 260 & 341 & 6 \\
11 & Constructing a Canal  bridge & 285 & 335 & 492 & 6 \\
12 & Construction works in the southern A1 lane & 113 & 128 & 189 & 9,10,11 \\
13 &Commissioning of the southern A1 new lane &
0 & 0 & 0 & 12 \\
14 & Producing parts of Railbridge part 1 & 223 & 251 & 366 & 6 \\
15 &  Producing parts of  Railbridge part 2 & 194 & 220 & 350 & 14 \\
16 &  Assembling a railway bridge on location & 559 & 674 & 971 & 14 \\
17 & Moving Railway Bridge in place during Train Free Period & 0 & 0 & 0 & 16 \\
18 & Road works northern A1 lane & 109 & 130 & 191 & 17 \\
19 & Commissioning of the Northern A1 new lane & 0 & 0 & 0 & 18 \\
20 & Road and construction works new part junction Diemen & 477 & 530 & 848 & 13 \\
21 & Build new viaducts A6 & 304 & 400 & 532 & 6 \\
22 & Build second Hollandse bridge & 286 & 340 & 459 & 6 \\
23 & Road and construction works junction Muiderberg & 716 & 930 & 1237 & 6 \\
24 & Road works eastern part A6 & 90 & 100 & 130 & 21,22 \\
25 & A6 East ready & 0 & 0 & 0 & 24 \\
26 & Reconstruction western part A6 & 324 & 400 & 532 & 25 \\
27 & Commissioning A6 & 0 & 0 & 0 & 18,23,26 \\
28 & Road works existing part Diemen junction & 71 & 90 & 130 & 19,20 \\
29 & Request Availability Certificate & 0 & 0 & 0 & 28 \\
30 & Assess and obtain Availability Certificate & 16 & 20 & 27 & 29 \\
31 & Demolition old A1 (part 1) & 54 & 61 & 91 & 30 \\
32 & Demolition old A1 (part 2) & 23 & 30 & 48 & 31 \\
33 & Greenery for old A1 & 77 & 90 & 129 & 31 \\
34 & Applications and obtaining partial completion certificates & 104 & 120 & 161 & 30 \\
35 & Request Completion Certificate & 0 & 0 & 0 & 33 \\
36 & Obtaining the Certificate of Completion & 17 & 20 & 27 & 35 \\
37 & Scheduled Completion Date & 0 & 0 & 0 & 36 \\
\hline
\end{tabular}}
\end{table}

\begin{table}[H]
    \centering
    \caption{Strategic control application: Shared uncertainty factors.}
    \label{tab:case3 correlation}
    \resizebox{\textwidth}{!}{
        \begin{tabular}{clrrrr}
\hline  &  & \multicolumn{3}{c}{ Shared uncertainty (days) } & \\
\hline ID & Shared uncertainty factor & $a_i$&  $m_i$ &  $b_i$ & Relations with activities\\
\hline 1 & Weather condition 1 & -45 & 0 & 72 & 10,11 \\
2 & Soil composition & -50 & 0 & 100 & $21,22,23,7$ \\
3 &Crew performance & -10 & 0 & 50 & 12,23 \\
4 &Soil composition & -45 & 0 & 110 & 20,26 \\
5 &Equipment availability 1 & -20 & 0 & 100 & 15,16 \\
6 &Site availability & -5 & 0 & 100 & 16,20 \\
7 &Procurement, fabrication or assembly & -1 & 0 & 55 & 7,20 \\
8 &Project control and management & -20 & 0 & 50 & 8,9 \\
9 &Design or documentation accuracy & -5 & 0 & 15 & 32,33 \\
10 &Owner-driven changes & 0 & 0 & 45 & 18,20 \\
11 &Issues with contractor & -20 & 0 & 50 & 3,4 \\
12 &Issues with supplier & -20 & 0 & 100 & 7,14 \\
13 &Equipment availability 2 & -80 & 0 & 90 & 7,16 \\
14 &Weather condition 2 & -140 & 0 & 100 & 7,23 \\
\hline
\end{tabular}}
\end{table}

\begin{table}[H]
    \centering
    \caption{Dynamic control application: Risk events.}
    \label{tab:saa_risks}
    \resizebox{\textwidth}{!}{
        \begin{tabular}{clrrrrr}
\hline \multirow[b]{2}{*}{ ID } & \multirow[b]{2}{*}{ Risk event description } & \multicolumn{3}{c}{ Risk duration (days) } & \multirow[b]{2}{*}{Relation} & \multirow[b]{2}{*}{$p_e \quad$} \\
\hline & &$a_i$&  $m_i$&  $b_i$& & \\
\hline 1 & Preliminary design rejection, including extra scope of works & 96 & 105 & 119 & 3 & 0.20 \\
2 & EDP audit failure & 13 & 14 & 15 & 4 & 0.05\\
3 & Condition deviates from plan & 63 & 70 & 78 & 7 & 0.15 \\
4 & Unexpected gas conducts & 35 & 35 & 41 & 8 & 0.20 \\
5 & Lower consolidation rate than calculated for & 34 & 35 & 40 & 9 & 0.10 \\
6 & Piling machines break down & 14 & 14 & 15 & 10 & 0.10 \\
7 & Late delivery of prefab elements & 19 & 21 & 25 & 11 & 0.20 \\
8 & Dynamic traffic management equipment/software not functioning & 20 & 21 & 22 & 12 & 0.25 \\
9 & Production equipment failure & 20 & 21 & 23 & 14 & 0.05\\
10 & Construction site subsides & 13 & 14 & 15 & 15 & 0.05 \\
11 & Ancillary equipment failure & 33 & 35 & 41 & 16 & 0.10 \\
12 & Dynamic traffic management equipment/software not functioning & 20 & 21 & 21 & 18 & 0.25 \\
13 & Discovery of polluted soil & 13 & 14 & 14 & 20 & 0.05 \\
14 & Concrete casting failure & 13 & 14 & 14 & 21 & 0.05\\
15 & Main pillar subsides & 65 & 70 & 71 & 22 & 0.02 \\
16 & Discovery of polluted soil & 25 & 28 & 32 & 23 & 0.05\\
17 & Insufficient quality of base layer & 39 & 42 & 47 & 26 & 0.02 \\
18 & Discovery of asphalt with too high PAK percentage & 13 & 14 & 17 & 28 & 0.05 \\
19 & Additional scope of work (miscellaneous) & 130 & 140 & 160 & 30 & 0.01 \\
\hline
\end{tabular}}
\end{table}

\begin{table}[H]
    \centering
    \caption{Dynamic control application: Corrective measures (variables).}
    \label{tab:saa_measures}
    \resizebox{\textwidth}{!}{
        \begin{tabular}{clrrrrrrrrrr}
\hline
& & \multicolumn{3}{c}{Capacity [\textit{days}]} & \multicolumn{3}{c}{Cost [\texteuro] $(\eta=0.5)$} & \multicolumn{3}{c}{Nuisance} & \\
\hline
 ID & Mitigation description &  $a_i$&  $m_i$ &  $b_i$&  $a_i$&  $m_i$&  $b_i$&  $a_i$&  $m_i$&  $b_i$& Relation\\
\hline
$x_1$ & Extra engineering design office personnel & 99 & 101 & 101 & 118k & 120k & 120k & 0.00 & 0.00 & 0.00 & 3 \\
$x_2$ & Extra software design capacity & 14 & 14 & 14 & 30k & 30k & 30k & 0.00 & 0.00 & 0.00 & 4 \\
$x_3$ & Extra maintenance engineers & 103 & 127 & 127 & 136k & 150k & 150k & 0.91 & 1.00 & 1.00 & 7 \\
$x_4$ & Extra administrators for permitting & 43 & 51 & 57 & 44k & 48k & 50k & 0.00 & 0.00 & 0.00 & 8 \\
$x_5$ & Applying extra preloading material & 41 & 51 & 51 & 677k & 750k & 750k & 5.42 & 6.00 & 6.00 & 9 \\
$x_6$ & Adding extra onsite construction flow & 92 & 101 & 107 & 190k & 200k & 205k & 9.52 & 10.00 & 10.00 & 10 \\
$x_7$ & Extra prefab construction capacity & 117 & 127 & 129 & 144k & 150k & 151k & 0.96 & 1.00 & 1.01 & 11 \\
$x_8$ & Extra M\&E engineers & 51 & 51 & 51 & 60k & 60k & 60k & 3.00 & 3.00 & 3.00 & 12 \\
$x_9$ & Extra welding equipment and personnel & 53 & 64 & 64 & 90k & 100k & 100k & 5.45 & 6.00 & 6.00 & 14 \\
$x_{10}$ & Extra temporary soil measures & 45 & 51 & 53 & 235k & 250k & 254k & 5.65 & 6.00 & 6.12 & 15 \\
$x_{11}$ & Ancillary on standby & 201 & 203 & 222 & 199k & 200k & 209k & 4.98 & 5.00 & 5.24 & 16 \\
$x_{12}$ & Extra M\&E engineers & 14 & 14 & 14 & 30k & 30k & 30k & 3.00 & 3.00 & 3.00 & 18 \\
$x_{13}$ & Extra excavation capacity & 96 & 101 & 101 & 121k & 125k & 125k & 9.71 & 10.00 & 10.00 & 20 \\
$x_{14}$ & Extra concrete workers/carpenters & 82 & 101 & 107 & 67k & 75k & 77k & 5.42 & 6.00 & 6.17 & 21 \\
$x_{15}$ & Temporary ancillary construction and rework & 70 & 76 & 84 & 1,442k & 1,500k & 1,576k & 7.69 & 8.00 & 8.41 & 22 \\
$x_{16}$ & Extra excavation capacity & 60 & 76 & 82 & 134k & 150k & 155k & 7.18 & 8.00 & 8.31 & 23 \\
$x_{17}$ & Extra asphalt equipment and personnel & 101 & 101 & 107 & 200k & 200k & 205k & 8.00 & 8.00 & 8.23 & 26 \\
$x_{18}$ & Extra removal works & 43 & 51 & 53 & 69k & 75k & 76k & 9.23 & 10.00 & 10.00 & 28 \\
$x_{19}$ & Extra equipment and personnel & 10 & 14 & 16 & 214k & 250k & 267k & 6.86 & 8.00 & 8.57 & 30 \\
$x_{20}$ & Applying extra preloading material during the night & 41 & 51 & 51 & 1,084k & 1,200k & 1,200k & 0.90 & 1.00 & 1.00 & 9 \\
$x_{21}$ & Automation of M\&E workflows & 30 & 30 & 30 & 120k & 120k & 120k & 0.00 & 0.00 & 0.00 & 12 \\
$x_{22}$ & Welding operation during the night & 45 & 53 & 53 & 166,k & 180k & 180k & 1.85 & 2.00 & 2.00 & 14 \\
$x_{23}$ & Extra temporary soil measures during the night & 32 & 40 & 42 & 315k & 350k & 358k & 0.90 & 1.00 & 1.03 & 15 \\
$x_{24}$ & Night shifts for concrete work & 67 & 80 & 91 & 91k & 100k & 106k & 3.68 & 4.00 & 4.28 & 21 \\
$x_{25}$ & Excavation during the night & 40 & 56 & 62 & 257k & 300k & 316k & 0.86 & 1.00 & 1.05 & 23 \\
$x_{26}$ & Removal works during the night & 33 & 41 & 43 & 81k & 90k & 92k & 0.90 & 1.00 & 1.02 & 28 \\
$x_{27}$ & Lane by lane asphalt work & 30 & 30 & 30 & 430k & 430k & 430k & 3.00 & 3.00 & 3.00 & 26 \\
\hline
\end{tabular}}
\end{table}


\bibliographystyle{elsarticle-harv} 
\bibliography{references}

\end{document}